\documentclass{amsart}

\usepackage{amsfonts}
\usepackage{txfonts}
\usepackage{amssymb}
\usepackage{enumitem}
\usepackage{graphicx}
\graphicspath{{figures/}}
\usepackage{epsfig,mathrsfs}
\usepackage[bf,SL,BF]{subfigure}
\usepackage{longtable}
\usepackage{booktabs}
\usepackage{cite}
\usepackage{diagbox}
\usepackage{caption}
\usepackage{float}
\usepackage{epstopdf}
\usepackage{latexsym, bm}
\usepackage{bbm}
\usepackage{algorithm, algorithmic}
\usepackage{cleveref}

\usepackage{amsthm}

\newcommand{\Rb}{\mathbb{R}}
\newcommand{\Cb}{\mathbb{C}}
\newcommand{\Nb}{\mathbb{N}}
\newcommand{\supp}{\operatorname{supp}}
\newcommand{\dist}{\operatorname{dist}}
\newcommand{\dd}{\,\mathrm{d}}
\newcommand{\X}{\mathbb{X}}

\newcommand{\Dcal}{\mathcal{D}}

\newtheorem{theorem}{Theorem}[section]
\newtheorem{lemma}[theorem]{Lemma}
\newtheorem{proposition}[theorem]{Proposition}

\theoremstyle{definition}
\newtheorem{definition}[theorem]{Definition}
\newtheorem{assumption}[theorem]{Assumption}
\theoremstyle{remark}
\newtheorem{remark}[theorem]{Remark}

\begin{document}

\title[Fractional Magnetic Schr\"odinger Equation with Drift]{%
  Simultaneous Inversion of Drift and Potential in the\\
  Fractional Magnetic Schr\"odinger Equation}

\author[Yun-Long Liu, Jin Wen, Jun-Liang Fu, and Xue-Juan Ren]{Yun-Long Liu$^{1}$, Jin Wen$^{1,2*}$, Jun-Liang Fu$^{1}$, and Xue-Juan Ren$^{3}$\\ \small \it 1. Department of Mathematics, Northwest Normal University, Lanzhou 730070, P. R. China.\\
2. Gansu Provincial Research Center for Basic Disciplines of Mathematics and Statistics, Lanzhou 730070, P. R. China. \\
3. College of Information Engineering, Tarim University, Alar 843300,
P. R. China}
	\date{\today}
	\thanks{{*}Corresponding author. Email: wenj@nwnu.edu.cn. The work described in this article was supported by the NSF of China (12261082, 11326234, 12501593), NSF of Gansu Province (145RJZA099), Scientific research project of Higher School in Gansu Province (2014A-012), Project of NWNU-LKQN2020-08, Graduate First Class Course Construction Project of Northwest Normal University in 2024 (2024KCLX01012),  Graduate Education and Curriculum Reform Project of Northwest Normal University in 2024 (2024KGLX01016), and Gansu Provincial Innovation Star Program for Graduate Students (2026CXZX-302).}

\begin{abstract}
We study the Calder\'on problem for the fractional magnetic Schr\"odinger equation with drift in the bivariate additive framework. When the magnetic potential $A$ is known and fixed, we prove that the drift $b$ and the scalar potential $q$ can be simultaneously and uniquely determined from exterior Dirichlet-to-Neumann measurements. When $A$ is unknown, we prove that the triple $(A,b,q)$ can only be recovered up to a natural gauge equivalence class, and that the drift $b$ is a gauge invariant, hence uniquely determined. The proofs rely on the Alessandrini identity, Runge approximation, and the unique continuation property for fractional magnetic operators. Furthermore, we have developed a numerical reconstruction framework that combines finite difference discretization, radial basis function surrogate modeling, and smoothed TV regularization. One-dimensional and two-dimensional numerical examples have verified the effectiveness of the proposed method and supported the theoretical uniqueness result.

\end{abstract}

\maketitle

\section{Introduction}\label{sec:intro}
Calder\'on proposed the classical Calder\'on inverse problem in 1980 \cite{Cal-1980}. Subsequently, this problem was extended to the non-local Calder\'on problem. In \cite{GSU-2016}, T. Ghosh et al. proposed for the first time to study the Calder\'on problem in the fractional Segal equation. Within two distinct theoretical scenarios, they established global uniqueness, showing that external boundary measurements are sufficient to uniquely identify the unknown potential coefficient. In \cite{CekicLinRuland2020}, Ceki\'c, Lin and R\"uland addressed the fractional Schr\"odinger equation augmented with a drift term. Their analysis yielded simultaneous unique reconstruction for both drift and scalar potential from infinitely many exterior measurements over a bounded spatial domain. Readers can consult \cite{Ghosh2020,BGU-2021,CR-2021,CMRU-2022,QU-2024,Fei-2024,BS-2024} for further progress concerning fractional Calder\'on type inverse problems.

Whereas the foregoing discussion concerns the nonlocal fractional Schr\"odinger setting, we now turn to the local magnetic Schr\"odinger equation (MSE), which may be viewed as a local magnetic counterpart.
Let $\Omega \subset \mathbb{R}^{n}$ stand for a bounded Lipschitz open domain.  Solutions to the Dirichlet boundary value problem for the MSE satisfy
\[
\begin{cases}
(-\Delta)_Au+qu:=-\Delta u-i\nabla\cdot(Au)-iA\cdot\nabla u+(|A|^2+q)u=0, &\quad \text{in}\ \Omega\\
u=f, &\quad \text{on}\ \partial\Omega,
\end{cases}
\]
where $f$ denotes prescribed boundary data, while $A$ and $q$ correspond respectively to vector magnetic potential and scalar potential inside the medium. Boundary measurement information is encoded within the DN map
\[
\Lambda_{A,q}: H^{1/2}(\partial \Omega)\to H^{-1/2}(\partial \Omega).
\]
Inverse problems associated with magnetic Schr\"odinger equations carry both theoretical mathematical significance and practical application value, since they generalize the original Calder\'on formulation \cite{Cal-1980}.  Originally motivated by mineral exploration, where it was used to locate valuable mineral resources, this problem has since found wide applications in electrical impedance tomography (EIT) \cite{Uhlmann-2009}, medical imaging \cite{Jossinet-1998}, cancer detection \cite{Zou+Guo-2003}, vital-sign monitoring \cite{Isaacson+Newell+Goble+Cheney-1990}, and other areas (see, e.g., \cite{Adler+Holder-2021}), as well as in engineering \cite{Hallaji-2014}. Configurations featuring non-vanishing magnetic fields are particularly noteworthy, given their tight connection with fixed-energy inverse scattering theory \cite{Nakamura-1995}. The MSE also arises naturally in various mathematical problems, including isotropic elasticity \cite{Nakamura+Uhlmann-1994,Nakamura+Uhlmann-2003}, special limiting forms of Maxwell and Schr\"odinger systems \cite{McDowall-2000,Eskin-2001}, Dirac equations \cite{Nakamura+Tsuchida-2000}, and the Stokes system \cite{Heck-2007}.

Research on fractional-order magnetic Schr\"odinger equations (FMSE) has also made some initial progress in recent years. In \cite{DAveniaSquassina2018}, d'Avenia and Squassina studied variational minimization problems involving non-local operators coupled with external magnetic potentials, establishing solution existence results via concentration compactness arguments. In \cite{Covi2020}, the magnetic potential $A$ is a binary function $A(x,y)$ defined on $\mathbb{R}^{n}\times\mathbb{R}^{n}$. The fractional magnetic gradient is given by $\nabla_A^s u(x,y)=\nabla^s u(x,y)+A(x,y)u(x,y)$. Global uniqueness results were derived: unknown electromagnetic coefficients inside a bounded domain can be fixed up to natural gauge freedom given infinitely many solution measurements collected over arbitrary exterior open subsets. The proof technique builds upon the Alessandrini identity together with Runge approximation theory, extending known findings for standard fractional Laplacian models.

In \cite{Covi+Monkkonen+Railo-2021},  Covi et al. generalized this bivariate formulation towards higher-order fractional magnetic Schr\"odinger systems for fractional powers $s\in\mathbb R^+\setminus\mathbb Z$, the magnetic potential being a binary tensor field of order $\lfloor s\rfloor + 1$, and proved the uniqueness of the higher-order fractional magnetic Schr\"odinger equations under gauge transformations. They also obtained the Runge approximation property and the fractional Poincar\'e inequality. This confirms that the binary framework is not limited to $s\in(0,1)$, but naturally extends to higher-order non-local magnetic operators. In \cite{Li2020}, Li introduced fractional magnetic operators equipped with single-variable magnetic potential $A(x)$ and scalar electric potential $q(x)$. Making use of Runge approximation arguments, scalar potential $q$ can be reconstructed using partial exterior DN observations.

In this paper, we add a drift term $b\cdot\nabla u$ to the FMSE and study the inverse problem of recovering the coefficients from exterior DN observational data.
\begin{equation}\label{eq11}
\begin{cases}
((-\Delta)_A^s + b\cdot\nabla + q)u = 0 & \text{in } \Omega,\\
u = f &\text{in } \Omega_e.
\end{cases}
\end{equation}
where \((-\Delta)_A^s\) is the magnetic fractional Laplacian in the bivariate framework, defined in Section~\ref{sec:prelim}. The domain $\Omega\subset\mathbb R^n$ stands for a non empty bounded connected open set with fractional exponent $0<s<1$. Coefficient spaces read $b\in W^{1-s,\infty}(\Omega)^n$, $q\in L^\infty(\Omega)$, and exterior trace data satisfies $f\in L^2(\Omega)$. The fractional Laplacian $(-\Delta)^s$ is defined via the hypersingular integral
\begin{equation}
\label{eq12}
(-\Delta)^s u(x) = c_{n,s} \, \mathrm{P.V.} \int_{\mathbb{R}^n} \frac{u(x) - u(\mathbf{y})}{|x - \mathbf{y}|^{n+2s}} \, d\mathbf{y},
\end{equation}
with the normalization constant
\begin{equation}\label{eq13}
c_{n,s} = \frac{2^{2s} s \, \Gamma(\frac{n}{2} + s)}{\pi^{n/2} \Gamma(1 - s)}.
\end{equation}
Owing to the non-local character of this operator, boundary conditions are imposed over the whole exterior domain $\Omega_e=\mathbb R^n\setminus\Omega$. Well posedness for the corresponding direct problem requires the following spectral assumption.
\begin{assumption}
\label{ass1}
If \(u\in H^{s_m}(\mathbb{R}^{n})\) solves
\begin{equation*}
\begin{cases}
((-\Delta)_A^{s} + b(x)\cdot\nabla + q(x))u(x) = 0, & x \in \Omega, \\
u(x) = 0, & x \in \Omega_e,
\end{cases}
\end{equation*}
Then we must have $u\equiv 0$. Equivalently, zero cannot be a Dirichlet eigenvalue associated with the operator $(-\Delta)_A^{s} + b\cdot\nabla + q$.
\end{assumption}
Under Assumption~\ref{ass1}, one can define the DN mapping $\Lambda^s_{A,b,q}: H^s(\Omega_e)\rightarrow(H^s(\Omega_e))^*$ relying on bilinear forms attached to the governing equation. We examine two distinct inverse reconstruction scenarios: first, magnetic potential $A$ is known a-priori, and we aim to retrieve drift $b$ and scalar potential $q$ inside $\Omega$ from $\Lambda^s_{A,b,q}$; second, all three coefficients $(A,b,q)$ remain unknown and need to be inferred from DN measurements. Due to intrinsic gauge equivalence $\sim$ present within FMSE models, full unique determination for triple $(A,b,q)$ is unattainable. Nevertheless, one may uniquely recover the complete gauge equivalence class of these potentials.
\begin{theorem}\label{thm:knownA}
Let \(\Omega\subset\Rb^n\), \(n\ge 1\), be a bounded Lipschitz domain,
\(s\in(1/2,1)\), and let \(A\in\mathcal{D}\) be fixed.
Let \(b_j\in W^{1-s,\infty}(\Omega)^n\) and \(q_j\in L^\infty(\Omega)\), \(j=1,2\),
satisfy the eigenvalue condition. If there exist non-empty open subsets
$W_1,W_2\subset\Omega_e$ such that
\[
  \Lambda_{A,b_1,q_1}^s[f]\big|_{W_2}
  =
  \Lambda_{A,b_2,q_2}^s[f]\big|_{W_2},
  \qquad \forall f\in C_c^\infty(W_1).
\]
Then \(b_1=b_2\) and \(q_1=q_2\) in \(\Omega\).
\end{theorem}

\begin{theorem}\label{thm:unknownA}
Let $\Omega\subset\mathbb R^n$, $n\ge 1$, denote a bounded Lipschitz domain with
$s\in(1/2,1)$. For triples $(A_j,b_j,q_j)$, $j=1,2$, where $A_j\in\mathcal{D}$,
$b_j\in W^{1-s,\infty}(\Omega)^n$, and $q_j\in L^\infty(\Omega)$, suppose the eigenvalue condition holds. Assume non empty open sets $W_1,W_2\subset\Omega_e$ satisfy that DN maps corresponding to
\[
  (-\Delta)_{A_j}^s u_j + b_j\cdot\nabla u_j + q_j u_j = 0
  \quad\text{in }\Omega, \qquad j=1,2,
\]
obey
\[
  \Lambda_{A_1,b_1,q_1}^s[f]\big|_{W_2}
  =
  \Lambda_{A_2,b_2,q_2}^s[f]\big|_{W_2},
  \qquad \forall f\in C_c^\infty(W_1).
\]
Then
\[
  (A_1,b_1,q_1)\sim(A_2,b_2,q_2),
\]
meaning the two triples match modulo gauge equivalence relation
$\sim$ whose formal definition appears in Section \ref{sec:gauge}. Equivalently,
\[
  (A_1)_{a\parallel}=(A_2)_{a\parallel},
  \qquad
  b_1=b_2,
  \qquad
  Q_1=Q_2,
\]
where
\[
  Q_j = q_j + \int_{\Rb^n}|A_j|^2\dd y + (\nabla\cdot)^s (A_j)_{s\parallel},
  \qquad j=1,2.
\]
In particular, the drift \(b\) is uniquely determined.
\end{theorem}
Our theoretical arguments build upon three foundational ingredients: integral identities derived for DN operators, the weak unique continuation property (WUCP), together with the Runge approximation property (RAP) stated next.
\begin{theorem}\label{thm:runge}
Let $\Omega\subset\mathbb R^n$, where $n\ge 1$is a bounded Lipschitz field, and $s\in(1/2,1)$. Let \(A\in\mathcal{D}\), \(b\in W^{1-s,\infty}(\Omega)^n\),
\(q\in L^\infty(\Omega)\), and assume the eigenvalue condition.
Let \(W\subset\Omega_e\) be a non-empty open set, and let \(P_{A,b,q}\) denote
the Poisson operator associated with the equation
\[
  (-\Delta)_A^s u + b\cdot\nabla u + q u = 0 \quad\text{in }\Omega,
  \qquad u|_{\Omega_e}=f.
\]
Then the following assertions hold.

\begin{enumerate}
  \item[(a)] The set
  \[
    \Dcal
    :=
    \bigl\{ P_{A,b,q}f - f : f\in C_c^\infty(W) \bigr\}
  \]
  is dense in \(L^2(\Omega)\) and in \(\widetilde{H}^s(\Omega)\).

  \item[(b)] Further suppose that $\partial\Omega$ is of class $C^\infty$, $A\in C_c^\infty(\Omega^2)$,
  \(b\in C_c^\infty(\Omega)^n\), and \(q\in C_c^\infty(\Omega)\) with
  \(\supp(b),\supp(q)\Subset\Omega\).
  Let \(d(x)\) be a positive \(C^\infty\)-function on \(\overline{\Omega}\)
  such that \(d(x)=\dist(x,\partial\Omega)\) near \(\partial\Omega\).
  Then for every \(g\in C^\infty(\overline{\Omega})\), every \(k\in\Nb\),
  and every \(\epsilon>0\), there exists an exterior datum \(f_\epsilon\in C_c^\infty(W)\)
  such that the corresponding solution \(u_\epsilon=P_{A,b,q}f_\epsilon\) satisfies
  \[
    \bigl\| d^{-s}(x)\,u_\epsilon - g \bigr\|_{C^k(\Omega)} < \epsilon.
  \]
\end{enumerate}
\end{theorem}
Qualitative Runge approximation constitutes an essential machinery when handling inverse problems governed by fractional Schr\"odinger type operators, originally introduced within reference \cite{GSU-2016}. To establish such property, authors in \cite{GSU-2016} leveraged unique continuation results obtained via Carleman estimates for fractional Schr\"odinger equations from \cite{Ruland-2015}. Relevant prior studies addressing unique continuation phenomena for fractional elliptic operators include \cite{FF-2014, FF-2015, Seo-2015, R-2017, Ghosh2020, Yu-2017}.

\section{Preliminaries}\label{sec:prelim}

\subsection{Fractional Sobolev spaces}
Let $\Omega \subset \Rb^n$ ($n\ge 1$) denote a bounded Lipschitz domain and write exterior domain $\Omega_e = \mathbb R^n \setminus \overline{\Omega}$. For $s\in(0,1)$ we define the fractional Sobolev spaces
\[
  H^s(\Rb^n) = W^{s,2}(\Rb^n), \qquad
  \widetilde{H}^s(\Omega) = \overline{C_c^\infty(\Omega)}^{\,H^s(\Rb^n)},
\]
and the quotient space $\X = H^s(\Rb^n)/\widetilde{H}^s(\Omega)$. For open subset $W\subset\Omega_e$, we write $C_c^\infty(W)$ for smooth functions compactly supported in $W$. These definitions are standard; see \cite[Section 2]{Covi2020} and \cite[Section 2.1]{CekicLinRuland2020}.

\subsection{Bivariate magnetic potential}
\begin{definition}\label{def:bivariate}
Let $A$ is a function
\[
  A:\Rb^n\times\Rb^n \to \Cb^n, \qquad (x,y)\mapsto A(x,y).
\]
For such $A$, we define its symmetric and antisymmetric parts by
\[
  A_s(x,y) = \frac{A(x,y)+A(y,x)}{2},
  \qquad
  A_a(x,y) = A(x,y)-A_s(x,y),
\]
so that \(A=A_a+A_s\). We further perform parallel perpendicular splitting relative to vector $x, y$,
\[
  A_\parallel(x,y) = \frac{A(x,y)\cdot(y-x)}{|y-x|^2}(y-x),
  \qquad
  A_\perp(x,y) = A(x,y)-A_\parallel(x,y).
\]
Define two point-wise $L^2$ integral quantities
\[
\mathcal{I}_1 A(x) \coloneqq \left(\int_{\Rb^n} |A(y,x)|^2 \dd y\right)^{1/2}, \qquad
\mathcal{I}_2 A(x) \coloneqq \left(\int_{\Rb^n} |A(x,y)|^2 \dd y\right)^{1/2}.
\]
\end{definition}
Whenever $A \in C_c^\infty$, these integrals remain finite. Since product $A_a\cdot A_s$ yields an antisymmetric scalar function and orthogonality $A_\parallel \perp A_\perp$ holds, we arrive at
\begin{equation}\label{eq21}
\begin{aligned}
\|A\|_{L^2}^2 &= \|A_a + A_s\|_{L^2}^2 = \|A_a\|_{L^2}^2 + \|A_s\|_{L^2}^2 + 2\langle A_a, A_s\rangle \\
&= \|A_a\|_{L^2}^2 + \|A_s\|_{L^2}^2 + 2\int_{\Rb^{2n}} A_a \cdot A_s \dd x \dd y = \|A_a\|_{L^2}^2 + \|A_s\|_{L^2}^2,
\end{aligned}
\end{equation}
and similarly \(\|A\|_{L^2}^2 = \|A_\parallel\|_{L^2}^2 + \|A_\perp\|_{L^2}^2\).
\begin{lemma}\cite{Covi2020}\label{lem:extension}
The symmetric-antisymmetric as well as parallel-perpendicular decompositions formulated in Definition~\ref{def:bivariate} hold almost everywhere for elements $A \in L^2(\mathbb R^{2n})$.
\end{lemma}
Let $\Omega \subset \Rb^n$ be open and let $\Omega_e = \Rb^n \setminus \overline{\Omega}$.
Let $p := \max\{2, n/(2s)\}$.
The vector potential $A : \Rb^{2n} \to \Cb^n$ and the scalar potential $q : \Rb^n \to \Rb$
are assumed to satisfy the following properties:
\begin{itemize}
  \item[(p1)] $\mathcal{I}_1 A, \mathcal{I}_2 A \in L^{2p}(\Rb^n)$;
  \item[(p2)] $A_{s\parallel} \in H^{sp-s}(\Rb^{2n}, \Cb^n)$;
  \item[(p3)] $A_{a\parallel}(x,y) \cdot (y-x) \ge 0$ for every pair $x,y\in\mathbb R^n$;
  \item[(p4)] $q \in L^p(\Omega)$;
  \item[(p5)] $A \in L^2(\Rb^{2n})$ and $\supp(A) \subseteq \Omega^2$.
\end{itemize}
We define the sets of admissible potentials
\[
\mathcal{D}_0 := \{A : \text{(p1)--(p3) hold}\}, \qquad
\mathcal{D} := \{A : \text{(p1)--(p3) and (p5) hold}\},
\]
\[
\mathcal{F}_0 := \{(A,q) : \text{(p1)--(p4) hold}\}, \qquad
\mathcal{F} := \{(A,q) : \text{(p1)--(p5) hold}\}.
\]
\begin{lemma}\cite{Covi2020}\label{lem:extension2}
 Let $A \in L^2(R^{2n})\cap \mathcal{D}_0$. Then there exist positive symmetric distributions
\(\sigma_1,\sigma_2\in\mathcal D'(\Rb^{2n})\) with \(\sigma_j\ge1\) such that
\[
  (A_j)_{a\parallel}=\alpha(\sigma_j-1),\qquad j=1,2.
\]
Consequently,
\[
  (A_1)_{a\parallel}-(A_2)_{a\parallel}=\alpha(\sigma_1-\sigma_2).
\]
\end{lemma}
\subsection{Magnetic fractional operator}
We begin by introducing the auxiliary kernel
\[
  \alpha(x,y)
  :=
  \frac{C_{n,s}^{1/2}}{\sqrt{2}}
  \frac{y-x}{|y-x|^{n/2+s+1}},
  \qquad x,y\in\Rb^n,\ x\neq y,
\]
where \(C_{n,s}>0\) is the normalization constant appearing in the definition of
\((-\Delta)^s\).
The fractional gradient of \(u\) is then defined as
\[
  \nabla^s u(x,y) := \alpha(x,y)\,(u(x)-u(y)),
\]
which generates symmetric, parallel vector fields over $\mathbb R^{2n}$. Equivalent expression,
\[
  \nabla^s u(x,y)
  =
  \frac{C_{n,s}^{1/2}}{\sqrt{2}}
  \frac{u(x)-u(y)}{|x-y|^{n/2+s+1}}(y-x).
\]
\begin{lemma}[{\cite[Lemma 2.7]{Covi2020}}]\label{lem:fractional-gradient-bounded}
The fractional gradient can be continuously extended as mapping
\[
  \nabla^s : H^r(\Rb^n) \to \langle D_x + D_y\rangle^{r-s} L^2(\Rb^{2n}),
\]
and if \(r \le s\) then also
\[
  \nabla^s : H^r(\Rb^n) \to H^{r-s}(\Rb^{2n}).
\]
\end{lemma}
The magnetic fractional Laplacian is then
\[
  (-\Delta)_A^s := (\nabla\cdot)_A^s(\nabla_A^s),
\]
where $(\nabla\cdot)_A^s$ is the adjoint of $\nabla_A^s$. By \cite[Lemma 3.3]{Covi2020}, we have the decomposition
\begin{equation}\label{eq22}
  (-\Delta)_A^s u = (-\Delta)^s u + 2\int_{\Rb^n} (A_{a\parallel}\cdot\nabla^s u)\,\dd y + Qu,
\end{equation}
where
\[
  Q := q + \int_{\Rb^n} |A|^2\,\dd y + (\nabla\cdot)^s A_{s\parallel}.
\]

\subsection{Well-posedness and DN map}
This subsection addresses well-posedness properties for direct problem~\eqref{eq11} as well as its adjoint counterpart.

Let $\Omega \subset \mathbb R^n$ be non-empty bounded connected Lipschitz domain, take $s\in(1/2,1)$, $A\in\mathcal A_0$, $b\in W^{1-s,\infty}(\Omega)^n$, $q\in L^\infty(\Omega)$. Given source term $F\in(\widetilde{H}^s(\Omega))^*$ and exterior trace $f\in H^s(\mathbb R^n)$, consider Dirichlet boundary value problem
\begin{equation}\label{eq23}
\begin{cases}
((-\Delta)_A^s + b\cdot\nabla + q)u = F & \text{in } \Omega,\\[2pt]
u - f \in \widetilde{H}^s(\Omega).
\end{cases}
\end{equation}

Given arbitrary open non-empty domain $\mathcal O\subset\mathbb R^n$ and $u,v\in L^2(\mathcal O)$, inner-product notation
\[
(u,v)_{\mathcal{O}} := \int_{\mathcal{O}} u v \dd x.
\]
For test functions $u,v\in C_c^\infty(\mathbb R^n)$, define main bilinear form
\begin{equation}\label{eq24}
B_{A,b,q}^s(u,v) := B_{A,q}^s(u,v) +  \langle b\cdot\nabla u,v\rangle_{\Omega}.
\end{equation}
where
\begin{equation}\label{eq25}
  B_{A,q}^s(u,v) = \langle \nabla_A^s u,\nabla_A^s v\rangle_{\Rb^{2n}} + \langle qu,v\rangle_{\Omega}.
\end{equation}
Since bilinear form $B_{A,b,q}^s(\cdot,\cdot)$ lacks symmetry, we further introduce adjoint bilinear form
\begin{equation}\label{eq26}
\begin{aligned}
B_{A,b,q}^{s,*}(u^*,v^*) :&=  \langle \nabla_A^s u^*,\nabla_A^s v^*\rangle_{\Rb^{2n}} + \langle bu^*,\nabla v^*\rangle_{\Omega}+  \langle q u^*, v^*\rangle_{\Omega}\\
&=B_{A,q}^{s,*}(u^*,v^*) + \langle bu^*,\nabla v^*\rangle_{\Omega}.
\end{aligned}
\end{equation}
for $v^*, w^* \in C_c^\infty(\Rb^n)$, associated with adjoint Dirichlet system
\begin{equation}\label{eq27}
\begin{cases}
(-\Delta)_A^s u^* - \nabla\cdot(b u^*) + qu^* = F \quad \text{in } \Omega,\\
u^* - f^* \in \widetilde{H}^s(\Omega).
\end{cases}
\end{equation}

\begin{remark}\label{rem:discrepancy}
As in \cite[Remark 2.3]{CekicLinRuland2020}, there exists a minor mismatch between adjoint bilinear form~\eqref{eq26} and the weak formulation for adjoint PDE~\eqref{eq27}. Indeed, when testing against an exterior extension $v^* \notin \widetilde{H}^s(\Omega)$, the two differ by the boundary term
\[
\int_{\partial\Omega} (b\cdot\nu) u^* v^* \dd S.
\]
This boundary contribution is absent when $v^* \in \widetilde{H}^s(\Omega)$, and it does not arise from the magnetic part $(-\Delta)_A^s$, which is self-adjoint. Operator $\Lambda_{A,b,q}^{s,*}$ is built purely from adjoint bilinear form without incorporating such boundary integral.
\end{remark}

\begin{lemma}\cite{CekicLinRuland2020}\label{lem:poincare}
Let integer $n\ge1$ and $s\in(0,1)$. For bounded open domain $\Omega\subset\mathbb R^n$, there exists constant $C>0$ such that
\[
\|v\|_{L^2(\Omega)} \le C (\operatorname{diam}(\Omega))^s \|(-\Delta)^{s/2}v\|_{L^2(\mathbb R^n)}
\quad\text{for all } v\in \widetilde{H}^s(\Omega).
\]
\end{lemma}

We next establish boundedness property for bilinear form corresponding to fractional magnetic Schr\"odinger operator with drift.
\begin{lemma}\label{lem:boundedness}
Let $\Omega\subset\mathbb R^n$ denote bounded Lipschitz domain, $s\in(1/2,1)$, $A\in\mathcal A_0$, $b\in W^{1-s,\infty}(\Omega)^n$, $q\in L^\infty(\Omega)$. Then bilinear form
\[
B_{A,b,q}^s(u,v)
:= B_{A,q}^s(u,v) + \int_\Omega (b\cdot\nabla u)v\, dx
\]
admits continuous extension onto whole product space $H^s(\mathbb R^n)\times H^s(\mathbb R^n)$. That is, there is constant $C>0$ satisfying
\[
|B_{A,b,q}^s(u,v)| \le C \|u\|_{H^s(\mathbb R^n)} \|v\|_{H^s(\mathbb R^n)}
\]
for arbitrary $u,v\in H^s(\mathbb R^n)$.
\end{lemma}

\begin{proof}
First, we directly have
\begin{equation*}
\begin{aligned}
\bigl|B_{A,q}^s(u,v)\bigr|& = |(\nabla_A^s u,\nabla_A^s v) + \langle qu,v\rangle|\\
&\leq |(\nabla_A^s u,\nabla_A^s v)| + |\langle qu,v\rangle|\\
&\leq \|\nabla_A^s u\|_{L^2}\|\nabla_A^s v\|_{L^2} + \|qu\|_{H^{-s}}\|v\|_{H^s}\\
&\leq C\|u\|_{H^s}\|v\|_{H^s} + C\|q\|_{L^p}\|u\|_{H^s}\|v\|_{H^s}\\
&\leq C\|u\|_{H^s}\|v\|_{H^s},
\end{aligned}
\end{equation*}
$C>0$ is constant and independent of $u,v$. Next, we will analyze the drift term.
As $1-s\in(0,\frac12)$, we have $\widetilde{H}^{1-s}(\Omega)=H^{1-s}(\Omega)$ and $\|v\|_{\widetilde{H}^{1-s}(\Omega)}\le C\|v\|_{H^{1-s}(\Rb^n)}$.
Next, we choose $B\in W^{1-s,\infty}(\Rb^n)$ such that
\[
B|_{\Omega}=b,\quad \|B\|_{W^{1-s,\infty}(\Rb^n)}\le 2\|b\|_{W^{1-s,\infty}(\Omega)},
\]
and estimate
\begin{equation*}
\begin{aligned}
\biggl|\int_{\Omega} v(b\cdot\nabla u)\,\mathrm{d}x\biggr|
&\le \|b v\|_{\widetilde{H}^{1-s}(\Omega)}\|\nabla u\|_{H^{s-1}(\Omega)}\\
&\le C\|B v\|_{H^{1-s}(\Rb^n)}\|\nabla v\|_{H^{s-1}(\Rb^n)}\\
&\le C\|v\|_{H^{1-s}(\Rb^n)}\|b\|_{W^{1-s,\infty}(\Omega)}\|\nabla u\|_{H^{s-1}(\Rb^n)}\\
&\le C\|v\|_{H^{1-s}(\Rb^n)}\|b\|_{W^{1-s,\infty}(\Omega)}\|u\|_{H^{s}(\Rb^n)}\\
&\le C\|v\|_{H^{s}(\Rb^n)}\|b\|_{W^{1-s,\infty}(\Omega)}\|u\|_{H^{s}(\Rb^n)}.
\end{aligned}
\end{equation*}
By Kato-Ponce inequality, we have
\begin{equation*}
\begin{aligned}
\|B v\|_{H^{1-s}(\Rb^n)}
&\le C\bigl\| J^{1-s}(B v)\bigr\|_{L^2(\Rb^n)}\\
&\le C\Bigl(\|B\|_{L^\infty(\Rb^n)}\bigl\| J^{1-s}v\bigr\|_{L^2(\Rb^n)}
+\bigl\| J^{1-s}B\bigr\|_{L^\infty(\Rb^n)}\|v\|_{L^2(\Rb^n)}\Bigr)\\
&\le C\|B\|_{W^{1-s,\infty}(\Rb^n)}\|v\|_{H^{1-s}(\Rb^n)}\\
&\le 2C\|b\|_{W^{1-s,\infty}(\Omega)}\|v\|_{H^{1-s}(\Rb^n)}\\
&\le 2C\|b\|_{W^{1-s,\infty}(\Omega)}\|v\|_{H^{s}(\Rb^n)},
\end{aligned}
\end{equation*}
where $J^{1-s}:=(\Delta-1)^{\frac{1-s}{2}}$ and $0<1-s<\frac12$. Finally, based on all the estimations, we can obtain
\[
\bigl|B_{A,b,q}^s(u,v)\bigr|\le C\|u\|_{H^{s}(\Rb^n)}\|v\|_{H^{s}(\Rb^n)}.
\]
We can arrive at the boundedness associated with the bilinear form through a similar argument,
\[
\bigl|B_{A,b,q}^{s,*}(u^{*},v^{*})\bigr|\le C\|u^{*}\|_{H^{s}(\Rb^n)}\|v^{*}\|_{H^{s}(\Rb^n)}.
\]
\end{proof}

By the above definitions and estimates, we have
\begin{equation}\label{eq29}
B_{A,b,q}^s(u,v)=B_{A,b,q}^{s,*}(v,u),\qquad \text{for any }u,v\in H^{s}(\Rb^n).
\end{equation}

We say that
\(u\in H^s(\Rb^n)\) with \(u-f\in\widetilde{H}^s(\Omega)\)
is a weak solution of
\begin{equation}\label{eq210}
\begin{cases}
((-\Delta)_A^s + b\cdot\nabla + q)u = F & \text{in } \Omega,\\[2pt]
u - f \in \widetilde{H}^s(\Omega),
\end{cases}
\end{equation}
if for all \(w\in\widetilde{H}^s(\Omega)\) we have
\[
B_{A,b,q}^s(u,w) = F(w).
\]
Weak solution of the adjoint problem
\begin{equation}\label{eq211}
\begin{cases}
(-\Delta)_A^s u^* - \nabla\cdot(b u^*) + q u^* = F^* & \text{in } \Omega,\\[2pt]
u^* - f^* \in \widetilde{H}^s(\Omega),
\end{cases}
\end{equation}
is defined analogously using \(B_{A,b,q}^{s,*}(\cdot,\cdot)\).

Next, we prove the well-posedness of the Dirichlet problem for the magnetic fractional Schr\"odinger equation with drift terms.

\begin{proposition}\label{prop:spectral-wellposed}
Let \(\Omega\subset\Rb^n\), \(n\ge 1\), \(s\in(1/2,1)\),
and let \(A\in\mathcal{D}_0\), \(b\in W^{1-s,\infty}(\Omega)^n\), \(q\in L^\infty(\Omega)\).
Let \(B_{A,b,q}^s(\cdot,\cdot)\) and \(B_{A,b,q}^{s,*}(\cdot,\cdot)\) be the bilinear and adjoint bilinear forms defined above. Then the following assertions hold.

\begin{enumerate}
  \item[(a)] There is a discrete set \(\Sigma=\{\lambda_j\}_{j=1}^\infty\subset\Rb\) with
  \(\lambda_1\le\lambda_2\le\cdots\to\infty\) such that for every \(\lambda\in\Rb\setminus\Sigma\),
  every \(f\in H^s(\Rb^n)\), and every \(F\in(\widetilde{H}^s(\Omega))^*\), there exists a unique
  \(u\in H^s(\Rb^n)\) satisfying
  \[
    B_{A,b,q}^s(u,w)-\lambda(u,w)_\Omega=F(w)
    \qquad\text{for all }w\in\widetilde{H}^s(\Omega),
  \]
  together with \(u-f\in\widetilde{H}^s(\Omega)\). Furthermore,
  \[
    \|u\|_{H^s(\Rb^n)}
    \le C\bigl(\|f\|_{H^s(\Rb^n)}+\|F\|_{(\widetilde{H}^s(\Omega))^*}\bigr),
  \]
  where \(C>0\) depends only on \(n\), \(s\), \(\|b\|_{W^{1-s,\infty}(\Omega)}\),
  \(\|q\|_{L^\infty(\Omega)}\), the admissible potential \(A\), and \(\lambda\).
  Moreover, \(\Sigma\) is bounded from below; more precisely,
  \[
    \Sigma\subset\bigl(-C(1+\|b\|_{W^{1-s,\infty}(\Omega)}+\|q\|_{L^\infty(\Omega)}
+\|A\|_{L^2(\Rb^{2n})}),\infty\bigr).
  \]

  \item[(b)] Similarly, for the bilinear form, there exists a discrete set
  \(\Sigma^*\subset\Rb\) such that for every \(\lambda^*\in\Rb\setminus\Sigma^*\) and every
  \(f^*\in H^s(\Rb^n)\), there is a unique \(u^*\in H^s(\Rb^n)\) with
  \[
    B_{A,b,q}^{s,*}(u^*,w^*)-\lambda^*(u^*,w^*)_\Omega=F^*(w^*)
    \qquad\text{for all }w^*\in\widetilde{H}^s(\Omega),
  \]
  and \(u^*-f^*\in\widetilde{H}^s(\Omega)\).
\end{enumerate}
\end{proposition}

\begin{proof}
We first prove (a). The argument relies on the spectral theorem for compact operators.
To begin, we may assume without loss of generality that $f=0$; indeed, setting
$v=u-f\in\widetilde{H}^s(\Omega)$ transfers the exterior datum into a modified functional
$\widetilde F\in(\widetilde{H}^s(\Omega))^*$. Thus it remains to solve for $v$. Since boundedness
of $B_{A,b,q}^s$ has already been established, it suffices to verify coercivity.

For $\psi\in C_c^\infty(\Rb^n)$ we recall the interpolation estimate
\[
  \|\psi\|_{H^{1-s}(\Rb^n)}
  \le C\|\psi\|_{L^2(\Rb^n)}^{\frac{2s-1}{s}}
  \|\psi\|_{H^s(\Rb^n)}^{\frac{1-s}{s}},
  \qquad \frac12<s<1.
\]
Using this and the Lipschitz property of $\Omega$, for $\psi,v\in C_c^\infty(\Rb^n)$ we estimate
the drift contribution as
\begin{align*}
  \Bigl|\int_{\Omega}\psi\,(\nabla v\cdot b)\dd x\Bigr|
  &\le \|b\psi\|_{\widetilde{H}^{1-s}(\Omega)}\|\nabla v\|_{H^{s-1}(\Omega)}\\
  &\le C\|b\psi\|_{H^{1-s}(\Rb^n)}\|\nabla v\|_{H^{s-1}(\Rb^n)}\\
  &\le C\|\psi\|_{H^{1-s}(\Rb^n)}\|b\|_{W^{1-s,\infty}(\Omega)}
      \|v\|_{H^s(\Rb^n)}\\
  &\le C\|v\|_{H^s(\Rb^n)}\|b\|_{W^{1-s,\infty}(\Omega)}
      \|\psi\|_{L^2(\Rb^n)}^{\frac{2s-1}{s}}
      \|\psi\|_{H^s(\Rb^n)}^{\frac{1-s}{s}}\\
  &\le C_\epsilon\|b\|_{W^{1-s,\infty}(\Omega)}\|v\|_{H^s(\Rb^n)}
      \|\psi\|_{L^2(\Rb^n)}
      +\epsilon\|v\|_{H^s(\Rb^n)}\|\psi\|_{H^s(\Rb^n)},
\end{align*}
where $\epsilon>0$ is arbitrary and we used Young's inequality in the last line. By density,
the same estimate extends to $v,\psi\in\widetilde{H}^s(\Omega)$.

Since $q\in L^\infty(\Omega)$, $A\in\mathcal{D}_0$, and $(\nabla\cdot)^s A_{s\parallel}\in L^p(\Omega)$ with $p = \max\{2, n/(2s)\} \ge n/(2s)$ (as defined in \cite{Covi2020}), we have $Q\in L^p(\Omega)$, where $Q = q + \int_{\mathbb{R}^n}|A|^2 dy + (\nabla\cdot)^s A_{s\parallel}$.

By H\"older's inequality,
\[
\Bigl|\int_\Omega Q v^2\,dx\Bigr|
\le \|Q\|_{L^p(\Omega)} \|v\|_{L^{2p'}(\Omega)}^2,
\qquad \frac{1}{p}+\frac{1}{p'}=1.
\]
By the fractional Gagliardo-Nirenberg interpolation inequality, there exists an interpolation exponent $\theta\in(0,1)$ such that
\[
\|v\|_{L^{2p'}(\Omega)} \le C \|v\|_{H^s(\mathbb{R}^n)}^\theta \|v\|_{L^2(\Omega)}^{1-\theta}.
\]
Squaring both sides and applying Young's inequality, for any $\epsilon>0$ there exists a constant $C_\epsilon>0$ depending on $\epsilon$ and $\|Q\|_{L^p(\Omega)}$ such that
\[
\|Q\|_{L^p(\Omega)} \|v\|_{L^{2p'}(\Omega)}^2
\le \epsilon\|v\|_{H^s(\mathbb{R}^n)}^2 + C_\epsilon \|v\|_{L^2(\Omega)}^2.
\]
By Lemma \ref{lem:poincare} and combining this with the lower bound $\langle\nabla^s v,\sigma\nabla^s v\rangle \ge \|(-\Delta)^{s/2}v\|_{L^2(\mathbb{R}^n)}^2$ and choosing $\epsilon>0$ sufficiently small, we obtain
\[
B_{A,q}^s(v,v)
\ge c\|v\|_{H^s(\mathbb{R}^n)}^2 - M_0\|v\|_{L^2(\Omega)}^2
\]
for constants $c>0$ and $M_0\ge 0$ depending only on $n$, $s$, $\Omega$, $\|Q\|_{L^p(\Omega)}$ and the admissible potential $A$.

Combining above estimates
\[
\begin{aligned}
  B_{A,b,q}^s(v,v)
  &= B_{A,q}^s(v,v) + \int_\Omega (b\cdot\nabla v)v \dd x \\
  &\ge \|(-\Delta)^{s/2}v\|_{L^2(\Rb^n)}^2
      - M_0\|v\|_{L^2(\Omega)}^2
      - \Bigl|\int_\Omega (b\cdot\nabla v)v \dd x\Bigr| \\
  &\ge \|(-\Delta)^{s/2}v\|_{L^2(\Rb^n)}^2
      - M_0\|v\|_{L^2(\Omega)}^2
      - \epsilon\|v\|_{H^s(\Rb^n)}^2
      - C_\epsilon\|b\|_{W^{1-s,\infty}(\Omega)}\|v\|_{L^2(\Omega)}^2.
\end{aligned}
\]
Using \(\|v\|_{H^s(\Rb^n)}^2 = \|v\|_{L^2(\Rb^n)}^2 + \|(-\Delta)^{s/2}v\|_{L^2(\Rb^n)}^2\), substituting this into the previous inequality, we get
\[
\begin{aligned}
  B_{A,b,q}^s(v,v)
  &\ge (1-\epsilon)\|v\|_{H^s(\Rb^n)}^2
      - \bigl(M_0 + C_\epsilon\|b\|_{W^{1-s,\infty}(\Omega)}+1\bigr)\|v\|_{L^2(\Omega)}^2 \\
  &\ge c\|v\|_{H^s(\Rb^n)}^2
      - \mu\|v\|_{L^2(\Omega)}^2.
\end{aligned}
\]
Thus, there exist constants \(c>0\) and \(\mu\ge 0\), depending only on
\(n,s,\Omega,\|b\|_{W^{1-s,\infty}(\Omega)},\|q\|_{L^\infty(\Omega)}\) and \(A\), such that
\[
  B_{A,b,q}^s(v,v)+\mu\|v\|_{L^2(\Omega)}^2
  \ge c\|v\|_{H^s(\Rb^n)}^2
  \qquad\text{for all }v\in\widetilde{H}^s(\Omega),
\]
As a consequence, $B_{A,b,q}^s(\cdot,\cdot)+\mu(\cdot,\cdot)_{L^2(\Omega)}$ defines a scalar product on
$\widetilde{H}^s(\Omega)$. By the Riesz representation theorem, for every
$\widetilde F\in(\widetilde{H}^s(\Omega))^*$ there exists a unique
$v\in\widetilde{H}^s(\Omega)$ such that
\[
  B_{A,b,q}^s(v,\psi)+\mu(v,\psi)_{L^2(\Omega)}=\widetilde F(\psi)
  \qquad\text{for all }\psi\in\widetilde{H}^s(\Omega).
\]
Writing $v=G_\mu(\widetilde F)$ defines a bounded linear operator
$G_\mu:(\widetilde{H}^s(\Omega))^*\to\widetilde{H}^s(\Omega)$. Since the embedding
$\widetilde{H}^s(\Omega)\hookrightarrow L^2(\Omega)$ is compact, $G_\mu$ is compact when viewed
as an operator on $L^2(\Omega)$.

The equation
\[
  B_{A,b,q}^s(v,\cdot)-\lambda(v,\cdot)_{L^2(\Omega)}=\widetilde F(\cdot)
  \quad\text{on }\widetilde{H}^s(\Omega)
\]
is equivalent to $v=G_\mu((\mu+\lambda)v+\widetilde F)$. The spectral theorem for compact
operators now yields a countable, decreasing sequence $\{\gamma_j\}_{j=1}^\infty\subset[0,\infty)$
with $\gamma_j\to0$ such that $G_\mu-\gamma\operatorname{Id}$ is invertible whenever
$\gamma\notin\{\gamma_j\}$. Setting $\gamma_j=(\lambda_j+\mu)^{-1}$ gives a countable,
increasing sequence $\lambda_j\to\infty$. From this, the proof of the solvability of $B_{A,b,q}^s(\cdot,\cdot)$ is obtained.

For (b), the adjoint case, the argument is analogous. One only needs to verify the boundedness
of the dual drift term
\[
  \Bigl|\int_\Omega b v^*\cdot\nabla\psi^*\dd x\Bigr|,
\]
where $v^*=u^*-f^*\in\widetilde{H}^s(\Omega)$ and $\psi^*\in\widetilde{H}^s(\Omega)$ is an
arbitrary test function. Using the Kato-Ponce inequality, Young's inequality, and the interpolation estimate
, we obtain for any $\epsilon>0$
\[
  \Bigl|\int_\Omega b v^*\cdot\nabla\psi^*\dd x\Bigr|
  \le C_\epsilon\|b\|_{W^{1-s,\infty}(\Omega)}\|\psi^*\|_{H^s(\Rb^n)}
      \|v^*\|_{L^2(\Omega)}
  +\epsilon\|\psi^*\|_{H^s(\Rb^n)}\|v^*\|_{H^s(\Rb^n)}.
\]
Prove that the rest of the similar to part (a). Thus, the proof is completed.
\end{proof}
Now, we consider the homogeneous Dirichlet problem
\begin{equation}\label{eq212}
\begin{cases}
((-\Delta)_A^s + b\cdot\nabla + q)u = 0 & \text{in } \Omega,\\[2pt]
u - f \in \widetilde{H}^s(\Omega),
\end{cases}
\end{equation}
for some exterior datum \(f \in H^s(\Omega_e)\).
By the eigenvalue condition, together with Proposition~\ref{prop:spectral-wellposed},
ensures equation \eqref{eq212} in $u \in H^s(\mathbb {R}^n)$ on have a unique solution.
We emphasize that this solution depends on \(f\) only through its equivalence class modulo
\(\widetilde{H}^s(\Omega)\): indeed, replacing \(f\) by \(f+\varphi\) with
\(\varphi \in \widetilde{H}^s(\Omega)\) leaves the exterior datum unchanged, so that
\(u_{f+\varphi}=u_f\) by uniqueness.

Then, we introduce the quotient space
\[
  \X = H^s(\Rb^n) / \widetilde{H}^s(\Omega),
\]
And \(f \in H^s(\Rb^n)\) equivalence class of \([f]\).
For a bounded Lipschitz domain \(\Omega\), one has
\(\X = H^s(\Omega_e)\).
Based on the well-posedness of \eqref{eq212}, we can now rigorously define the DN map as follows.

\begin{definition}[DN map]\label{def:DN}
Let \(\Omega \subset \Rb^n\), \(n \ge 1\), \(s \in (1/2,1)\),
and let \(A \in \mathcal{D}_0\), \(b \in W^{1-s,\infty}(\Omega)^n\), \(q \in L^\infty(\Omega)\).
Let \(B_{A,b,q}^s(\cdot,\cdot)\) be the bilinear form defined in \eqref{eq24}.
The DN operator associated with this equation is defined as
\[
\begin{cases}
(-\Delta)_A^s u + b\cdot\nabla u + q u = 0 & \text{in } \Omega,\\[2pt]
u = f & \text{in } \Omega_e,
\end{cases}
\]
as
\[
  \Lambda_{A,b,q}^s \colon \X \to \X^{*},
  \qquad
  \langle \Lambda_{A,b,q}^s[f],[g]\rangle = B_{A,b,q}^s(u_f,g),
\]
where \(f,g \in H^s(\Rb^n)\) and \(u_f\) denotes the unique weak solution to the above
equation with exterior datum \(f\).
\end{definition}

\begin{remark}\label{rem:DN-well-defined}
The definition above is independent of the chosen representatives. Indeed, for any
\(\varphi \in \widetilde{H}^s(\Omega)\) we have
\[
  B_{A,b,q}^s(u_f,g+\varphi) = B_{A,b,q}^s(u_f,g),
\]
since \(\varphi\) vanishes on \(\Omega_e\) and \(u_f\) solves the equation in \(\Omega\).
Similarly, for any \(\varphi \in \widetilde{H}^s(\Omega)\),
\[
  B_{A,b,q}^s(u_{f+\varphi},g) = B_{A,b,q}^s(u_f,g),
\]
because \(u_{f+\varphi}=u_f\) by uniqueness. Consequently, the pairing
\(\langle \Lambda_{A,b,q}^s[f],[g]\rangle\) is well-defined on \(\X \times \X\).
Moreover, the boundedness of \(\Lambda_{A,b,q}^s \colon \X \to \X^{*}\) follows from the
multiplier estimate and the well-posedness established in
Proposition~\ref{prop:spectral-wellposed}.
\end{remark}
\begin{remark}\label{rem:DN-formal}
A formal computation using Definition~\ref{def:DN} shows that
\begin{equation}\label{eq213}
\begin{aligned}
  \langle \Lambda_{A,b,q}^s[f],[g]\rangle
  &= B_{A,b,q}^s(u_f,g) \\
  &= \langle \nabla_A^s u_f,\nabla_A^s g\rangle_{\Rb^{2n}}
     + \langle q u_f, g\rangle_{\Omega}
     + \langle b\cdot\nabla u_f, g\rangle_{\Omega} \\
  &= \langle (-\Delta)_A^s u_f, g\rangle_{\Omega_e}
     + \bigl\langle (-\Delta)_A^s u_f + b\cdot\nabla u_f + q u_f, g\bigr\rangle_{\Omega} \\
  &= \int_{\Omega_e} g\,(-\Delta)_A^s u_f \dd x,
\end{aligned}
\end{equation}
where \(u_f\) is weak solution to \eqref{eq212} with exterior datum \(f\).
From \eqref{eq213} one may formally read off
\[
  \Lambda_{A,b,q}^s[f] = (-\Delta)_A^s u_f\bigr|_{\Omega_e},
\]
which is the nonlocal analogue of the classical DN map.
A rigorous justification of this identity requires additional regularity assumptions
on the coefficients and the data, we refer to \cite{GSU-2016} for the
corresponding details in the magnetic case.
\end{remark}

Next, we analyze the adjoint problem with vanishing source term,
\begin{equation}\label{eq214}
\begin{cases}
(-\Delta)_A^s u^* - \nabla\cdot(b u^*) + q u^* = 0 & \text{in } \Omega,\\[2pt]
u^* - f^* \in \widetilde{H}^s(\Omega).
\end{cases}
\end{equation}
By the eigenvalue condition and Proposition~\ref{prop:spectral-wellposed},
\eqref{eq214} admits a unique solution \(u^*\) for every
\(f^* \in H^s(\Rb^n)\). We define the DN map
\(\Lambda_{A,b,q}^{s,*}\) associated with the adjoint bilinear form
\(B_{A,b,q}^{s,*}(\cdot,\cdot)\) by
\[
  \Lambda_{A,b,q}^{s,*} \colon \X \to \X^{*},
  \qquad
  \langle \Lambda_{A,b,q}^{s,*}[f],[g]\rangle
  = B_{A,b,q}^{s,*}(u_f^*,g),
\]
where \(u_f^* \in H^s(\Rb^n)\) is the unique weak solution of \eqref{eq214}
with exterior datum \(f\).

Let \(\Lambda_{A,b,q}^{s,\prime} \colon \X \to \X^{*}\) denote the adjoint of
\(\Lambda_{A,b,q}^s\) with respect to the duality pairing on \(\X\), i.e.
\[
  \langle \Lambda_{A,b,q}^s f, g\rangle
  = \langle f, \Lambda_{A,b,q}^{s,\prime} g\rangle,
  \qquad f,g \in H^s(\Rb^n).
\]
The following result shows that the DN map \(\Lambda_{A,b,q}^{s,\prime}\)
coincides with the DN map \(\Lambda_{A,b,q}^{s,*}\) of the adjoint equation.

\begin{lemma}\label{lem:adjoint-DN}
Let \(\Omega \subset \Rb^n\), \(n \ge 1\), \(s \in (1/2,1)\), \(A \in \mathcal{D}_0\), \(b \in W^{1-s,\infty}(\Omega)^n\),
\(q \in L^\infty(\Omega)\), and assume the eigenvalue condition holds.
Let \(B_{A,b,q}^s(\cdot,\cdot)\), \(B_{A,b,q}^{s,*}(\cdot,\cdot)\) and
\(\Lambda_{A,b,q}^s\), \(\Lambda_{A,b,q}^{s,*}\) be defined as above.
Then, for all \(f,g \in \X\),
\begin{equation}\label{eq215}
  \langle \Lambda_{A,b,q}^s f, g\rangle
  = \langle f, \Lambda_{A,b,q}^{s,*} g\rangle.
\end{equation}
\end{lemma}
\begin{proof}
Let \(u_f\) solve \eqref{eq212} with exterior datum \(f \in \X\), and let
\(u_g^*\) solve the adjoint equation \eqref{eq214} with exterior
datum \(g \in \X\). By the duality relation between \(B_{A,b,q}^s\) and
\(B_{A,b,q}^{s,*}\), and by the definitions of \(\Lambda_{A,b,q}^s\) and
\(\Lambda_{A,b,q}^{s,*}\), we compute
\[
  \langle \Lambda_{A,b,q}^s f, g\rangle
  = B_{A,b,q}^s(u_f,u_g^*)
  = B_{A,b,q}^{s,*}(u_g^*,u_f)
  = \langle \Lambda_{A,b,q}^{s,*} g, f\rangle
  = \langle f, \Lambda_{A,b,q}^{s,*} g\rangle.
\]
\end{proof}

\begin{lemma}[Alessandrini identity]\label{lem:alessandrini}
Let \(\Omega \subset \Rb^n\), \(n \ge 1\), be a bounded Lipschitz domain,
\(s \in (1/2,1)\), and let \(A_j \in \mathcal{D}_0\), \(b_j \in W^{1-s,\infty}(\Omega)^n\),
\(q_j \in L^\infty(\Omega)\) for \(j=1,2\), satisfying the eigenvalue condition.
Let \(f_1, f_2 \in \mathbb{X}\), we have
\begin{equation}\label{eq:alessandrini}
\begin{aligned}
  \langle (\Lambda_{A_1,b_1,q_1}^s - \Lambda_{A_2,b_2,q_2}^s) f_1, f_2 \rangle
  &= 2\left\langle \int_{\Rb^n} ((A_1)_{a\parallel} - (A_2)_{a\parallel})
     \cdot \nabla^s u_1 \dd y,\; u_2^* \right\rangle \\
  &\quad + \langle (Q_1 - Q_2)u_1, u_2^* \rangle
     + \int_\Omega (b_1 - b_2)\cdot\nabla u_1\, u_2^* \dd x,
\end{aligned}
\end{equation}
where  \(u_1 \in H^s(\Rb^n)\) solve
\(
  (-\Delta)_{A_1}^s u_1 + b_1\cdot\nabla u_1 + q_1 u_1 = 0 \quad\text{in }\Omega,
\)
and \(u_2^* \in H^s(\Rb^n)\) solves the adjoint equation
\(
  (-\Delta)_{A_2}^s u_2^* - \nabla\cdot(b_2 u_2^*) + q_2 u_2^* = 0 \quad\text{in }\Omega,
\) with $u_1-f_1 \in \widetilde{H}^s(\Omega)$ and $u_2-f_2 \in \widetilde{H}^s(\Omega)$.
\end{lemma}

\begin{proof}
According to the above definition, we get
\begin{align*}
  &\langle (\Lambda_{A_1,b_1,q_1}^s - \Lambda_{A_2,b_2,q_2}^s) f_1, f_2 \rangle \\
  &= \langle \Lambda_{A_1,b_1,q_1}^s f_1, f_2 \rangle
     - \langle \Lambda_{A_2,b_2,q_2}^s f_1, f_2 \rangle \\
  &= B_{A_1,b_1,q_1}^s(u_1, f_2)
     - B_{A_2,b_2,q_2}^{s,*}(u_2^*, f_1) \\
  &= B_{A_1,b_1,q_1}^s(u_1, u_2^*)
     - B_{A_2,b_2,q_2}^{s,*}(u_2^*, u_1) \\
  &= B_{A_1,q_1}^s(u_1, u_2^*) + \int_\Omega (b_1\cdot\nabla u_1) u_2^* \dd x
     - B_{A_2,q_2}^s(u_2^*, u_1) - \int_\Omega b_2 u_2^* \cdot \nabla u_1 \dd x \\
  &= B_{A_1,q_1}^s(u_1, u_2^*) - B_{A_2,q_2}^s(u_1, u_2^*)
     + \int_\Omega (b_1 - b_2)\cdot\nabla u_1\, u_2^* \dd x \\
  &= 2\left\langle \int_{\Rb^n} ((A_1)_{a\parallel} - (A_2)_{a\parallel})
     \cdot \nabla^s u_1 \dd y,\; u_2^* \right\rangle
     + \langle (Q_1 - Q_2)u_1, u_2^* \rangle \\
  &\quad + \int_\Omega (b_1 - b_2)\cdot\nabla u_1\, u_2^* \dd x.
\end{align*}
\end{proof}

\section{Runge approximation}\label{Runge}
This section investigates Runge-type approximation behaviour for solutions of governing equation~\eqref{eq11}. We first introduce the weak version of unique continuation principle.
\begin{lemma}[Weak unique continuation]\label{lem:WUCP}
Let \((A,q)\in\mathcal P\) and \(b\in W^{1-s,\infty}(\Omega)^n\).
Assume \(W\subset\Omega_e\) is open and \(u\in H^s(\Rb^n)\) satisfies
\[
  u=0 \quad\text{and}\quad (-\Delta)_A^s u + b\cdot\nabla u + q u = 0 \quad\text{in }W.
\]
Then \(u=0\) in \(\Omega\).
\end{lemma}

\begin{proof}
Since \((A,q)\in\mathcal P\), we have \(\supp(A)\subseteq\Omega^2\), hence \(A_{a\parallel}=0\) on \(W\times\Rb^n\)
and \((-\Delta)_A^s u=(-\Delta)^s u\) in \(W\).
Moreover, \(u=0\) in \(W\) implies \(b\cdot\nabla u=0\) and \(q u=0\) in \(W\).
Thus the equation reduces to \((-\Delta)^s u=0\) in \(W\).
By virtue of the weak unique continuation property of the fractional Laplacian operator (refer to \cite[Theorem 1.2]{GSU-2016}), it follows that \(u = 0\) on \(\Omega\).
\end{proof}

\begin{lemma}[Runge approximation]\label{lem:runge}
Let \(\Omega\subset\Rb^n\), \(n\ge 1\), be a bounded Lipschitz domain and
\(s\in(1/2,1)\). Let \(A\in\mathcal{D}\), \(b\in W^{1-s,\infty}(\Omega)^n\),
\(q\in L^\infty(\Omega)\), and assume the eigenvalue condition.
Let \(W\subset\Omega_e\) be a non-empty open set, and let \(P_{A,b,q}\) denote
the Poisson operator associated with the equation
\[
  (-\Delta)_A^s u + b\cdot\nabla u + q u = 0 \quad\text{in }\Omega,
  \qquad u|_{\Omega_e}=f.
\]
Let
\[
  \mathcal{D}
  :=
  \bigl\{ P_{A,b,q}f - f : f\in C_c^\infty(W) \bigr\}.
\]
Then:
\begin{enumerate}
  \item[(a)] \(\mathcal{D}\) is dense in \(L^2(\Omega)\);
  \item[(b)] \(\mathcal{D}\) is dense in \(\widetilde{H}^s(\Omega)\).
\end{enumerate}
\end{lemma}

\begin{proof}
It suffices to prove (b), since \(\widetilde{H}^s(\Omega)\subset L^2(\Omega)\) with dense inclusion
and (a) then follows from (b) by density.

Applying Hahn-Banach separation principle, we only need to show, any linear functional $F\in(\widetilde{H}^s(\Omega))^*$ which annihilates all elements from set $\mathcal D$ must be identically zero.

Let \(\phi\in\widetilde{H}^s(\Omega)\) be the unique solution of the adjoint problem
\[
  (-\Delta)_A^s\phi - \nabla\cdot(b\phi) + q\phi = F \quad\text{in }\Omega,
  \qquad \phi=0 \quad\text{in }\Omega_e.
\]
In weak form, this reads
\[
  B_{A,b,q}^{s,*}(\phi,w)=F(w)
  \qquad\text{for all }w\in\widetilde{H}^s(\Omega).
\]
For \(f\in C_c^\infty(W)\), set \(u_f:=P_{A,b,q}f\in H^s(\Rb^n)\). Then
\(u_f-f\in\widetilde{H}^s(\Omega)\) and
\begin{align*}
  F(u_f-f)
  &= B_{A,b,q}^{s,*}(\phi,u_f-f)\\
  &= B_{A,b,q}^s(u_f-f,\phi)\\
  &= B_{A,b,q}^s(u_f,\phi)-B_{A,b,q}^s(f,\phi)\\
  &= -B_{A,b,q}^s(f,\phi),
\end{align*}
where we used the duality relation \(B_{A,b,q}^{s,*}(u,v)=B_{A,b,q}^s(v,u)\),
so \(B_{A,b,q}^s(u_f,\phi)=0\) for \(\phi\in\widetilde{H}^s(\Omega)\).
Since \(F\) vanishes on \(\mathcal{D}\), we obtain
\[
  B_{A,b,q}^s(f,\phi)=0
  \qquad\text{for all }f\in C_c^\infty(W).
\]
We now simplify this identity term by term. Recall that
\[
  B_{A,b,q}^s(f,\phi)
  =
  \langle \nabla_A^s f,\nabla_A^s\phi\rangle_{\Rb^{2n}}
  +
  \langle qf,\phi\rangle_{\Omega}
  +
  \langle b\cdot\nabla f,\phi\rangle_{\Omega},
\]
and that
\[
  \nabla_A^s f=\nabla^s f+A(x,y)f(x),
  \qquad
  \nabla_A^s\phi=\nabla^s\phi+A(x,y)\phi(x).
\]
Hence
\begin{equation}\label{eq:expansion-B}
\begin{aligned}
  \langle \nabla_A^s f,\nabla_A^s\phi\rangle_{\Rb^{2n}}
  &=
  \langle \nabla^s f,\nabla^s\phi\rangle_{\Rb^{2n}}
  +
  \langle \nabla^s f,A\phi\rangle_{\Rb^{2n}}
  +
  \langle Af,\nabla^s\phi\rangle_{\Rb^{2n}}
  +
  \langle Af,A\phi\rangle_{\Rb^{2n}}.
\end{aligned}
\end{equation}

Because \(f\in C_c^\infty(W)\), \(W\subset\Omega_e\) and  \(\phi\in\widetilde{H}^s(\Omega)\), so we have
$f=0,\ \nabla f=0\ \text{on }\Omega,\ \phi=0\ \text{on }\Omega_e.$\\
Since the integration in \(\langle qf,\phi\rangle_{\Omega}\) is restricted to \(\Omega\),
and \(f=0\) on \(\Omega\), we have
\[
  \langle qf,\phi\rangle_{\Omega}
  =
  \int_\Omega q(x)f(x)\phi(x)\dd x
  =
  0.
\]
Similarly, since \(\nabla f=0\) on \(\Omega\),
\[
  \langle b\cdot\nabla f,\phi\rangle_{\Omega}
  =
  \int_\Omega b(x)\cdot\nabla f(x)\,\phi(x)\dd x
  =
  0.
\]
Next, we analysis \eqref{eq:expansion-B},

For \(\langle Af,\nabla^s\phi\rangle_{\Rb^{2n}}\) and \(\langle Af,A\phi\rangle_{\Rb^{2n}}\), since \(\supp(A)\subseteq\Omega^2\) by (p5), the integrand can be nonzero only when
\(x,y\in\Omega\). But \(f=0\) on \(\Omega\). Therefore, both of these are zero.

For \(\langle \nabla^s f,A\phi\rangle_{\Rb^{2n}}\),
  \[
    \langle \nabla^s f,A\phi\rangle_{\Rb^{2n}}
    =
    \int_{\Rb^{2n}} \nabla^s f(x,y)\cdot A(x,y)\phi(x)\dd x\dd y.
  \]
  If \(A(x,y)\neq0\), then \((x,y)\in\Omega^2\) by (p5), so \(x,y\in\Omega\). But then
  \[
    \nabla^s f(x,y)
    =
    (f(x)-f(y))\,\alpha(x,y)
    =
    0,
  \]
  because \(f=0\) on \(\Omega\). Hence the integrand is zero and the integral vanishes.

Consequently, we get
\[
  B_{A,b,q}^s(f,\phi)
  =
  \langle \nabla^s f,\nabla^s\phi\rangle_{\Rb^{2n}}.
\]
By \((\nabla\cdot)^s\nabla^s=(-\Delta)^s\), we obtain
\[
  \langle \nabla^s f,\nabla^s\phi\rangle_{\Rb^{2n}}
  =
  \langle (\nabla\cdot)^s\nabla^s f,\phi\rangle_{\Rb^n}
  =
  \langle (-\Delta)^s f,\phi\rangle_{\Rb^n}
  =
  \langle f,(-\Delta)^s\phi\rangle_{\Rb^n}.
\]
Therefore, we have
\[
  \langle f,(-\Delta)^s\phi\rangle_{\Rb^n}=0
  \qquad\text{for all }f\in C_c^\infty(W).
\]
Since \(C_c^\infty(W)\) is dense in \(L^2(W)\), this implies
\[
  (-\Delta)^s\phi=0
  \quad\text{in }W.
\]

By \(\phi\in\widetilde{H}^s(\Omega)\), we have \(\phi=0\) in \(\Omega_e\), and in
particular \(\phi=0\) in \(W\). Moreover, in \(W\subset\Omega_e\) we have
\(A_{a\parallel}=0\) by \(\supp(A)\subseteq\Omega^2\), so that
\[
  (-\Delta)_A^s\phi = (-\Delta)^s\phi = 0
  \quad\text{in }W.
\]
Together with \(b\cdot\nabla\phi=0\) and \(q\phi=0\) in \(W\) (as \(\phi=0\) there),
we conclude
\[
  \phi=0
  \quad\text{and}\quad
  (-\Delta)_A^s\phi + b\cdot\nabla\phi + q\phi = 0
  \quad\text{in }W.
\]
By Lemma~\ref{lem:WUCP}, it follows that
\(\phi\equiv 0\) in \(\Omega\), and hence \(\phi\equiv 0\) in \(\Rb^n\). Therefore
\(F\equiv 0\), which concludes the proof.
\end{proof}
Theorem \ref{thm:runge}(a) has been proved. Theorem \ref{thm:runge}(b) is proved in the latter part of this paper.

\section{Uniqueness}\label{sec:uniquness}
In this section, we have proved Theorem \ref{thm:knownA} and Theorem \ref{thm:unknownA}.

\subsection{Gauge equivalences}\label{sec:gauge}
\(G\) are strictly for all functions \(\phi\in C^\infty(\Rb^n)\) Abelian group, makes \(\phi|_{\Omega_e}=1\).
For two triples \((A_1,b_1,q_1)\) and \((A_2,b_2,q_2)\), we introduce two equivalence relations.

\begin{definition}[Gauge \(\sim\)]\label{def:gauge-sim}
We say that \((A_1,b_1,q_1)\) and \((A_2,b_2,q_2)\) are \emph{\(\sim\)-gauge
equivalent}, written
\[
  (A_1,b_1,q_1)\sim(A_2,b_2,q_2),
\]
if
\begin{equation}\label{eq:gauge-sim}
  (-\Delta)_{A_1}^s u+b_1\cdot\nabla u+q_1u
  =
  (-\Delta)_{A_2}^s u+b_2\cdot\nabla u+q_2u
  \qquad\forall u\in H^s(\Rb^n).
\end{equation}
\end{definition}

\begin{definition}[Gauge \(\approx\)]\label{def:gauge-approx}
We say that \((A_1,b_1,q_1)\) and \((A_2,b_2,q_2)\) are \emph{\(\approx\)-gauge
equivalent}, written
\[
  (A_1,b_1,q_1)\approx(A_2,b_2,q_2),
\]
if there exists \(\phi\in G\) such that
\begin{equation}\label{eq:gauge-approx}
  (-\Delta)_{A_1}^s(u\phi)+b_1\cdot\nabla(u\phi)+q_1u\phi
  =
  \phi\bigl((-\Delta)_{A_2}^s u+b_2\cdot\nabla u+q_2u\bigr)
  \qquad\forall u\in H^s(\Rb^n).
\end{equation}
\end{definition}

Both \(\sim\) and \(\approx\) are equivalence relations on the set of admissible
triples. Since \(\phi\equiv1\in G\), we have
\[
  (A_1,b_1,q_1)\sim(A_2,b_2,q_2)
  \ \Longrightarrow\
  (A_1,b_1,q_1)\approx(A_2,b_2,q_2).
\]

\begin{remark}\label{rem:gauge-approx}
For the classical magnetic Schr\"odinger equation (\(s=1\)), the relevant
gauge is \(\approx\): the substitution \(u\mapsto\phi u\) with \(\phi\in G\)
produces a nontrivial family of equivalent potentials.
For the FMSE with \(b=0\), Covi \cite[Lemma 3.10]{Covi2020} proved that
\emph{no} nontrivial \(\approx\)-gauge exists: \(\approx\) reduces to \(\sim\).
When a drift term \(b\neq0\) is present, the same conclusion persists: expanding
\(b\cdot\nabla(u\phi)=\phi\,b\cdot\nabla u+u\,b\cdot\nabla\phi\) in
\eqref{eq:gauge-approx} and comparing the coefficients of \(\nabla u\), one
obtains \(b_1=b_2\), after which the argument of
\cite[Lemma 3.10]{Covi2020} applies verbatim. Hence the FMSE with drift
possesses the gauge \(\sim\) but not the gauge \(\approx\), in sharp contrast
to its local counterpart.
\end{remark}

\subsection{Case I: Known magnetic potential}\label{sec:knownA}
\begin{proof}[\textbf{Proof of Theorem \ref{thm:knownA}}]
Since \(A_1=A_2=A\), the Alessandrini identity \eqref{eq:alessandrini} reduces to
\begin{equation}\label{eq:knownA-identity}
  0
  =
  \langle (q_1-q_2)u_1,u_2^*\rangle_\Omega
  +
  \int_\Omega (b_1-b_2)\cdot\nabla u_1\, u_2^* \dd x,
\end{equation}
where \(u_1\) solves
\[
  (-\Delta)_A^s u_1+b_1\cdot\nabla u_1+q_1 u_1=0\quad\text{in }\Omega,
  \qquad u_1-f_1\in\widetilde{H}^s(\Omega),
\]
and \(u_2^*\) solves the adjoint equation
\[
  (-\Delta)_A^s u_2^*-\nabla\cdot(b_2 u_2^*)+q_2 u_2^*=0\quad\text{in }\Omega,
  \qquad u_2^*-f_2\in\widetilde{H}^s(\Omega).
\]
We recover \(q\) and \(b\) in two steps.

\emph{ Step 1: Recovery of \(q\).}
Let \(\psi_2\in C_c^\infty(\Omega)\) be arbitrary. Choose a cutoff function
\(\psi_1\in C_c^\infty(\Omega)\) with \(\psi_1=1\) on \(\operatorname{supp}(\psi_2)\).
By Lemma~\ref{lem:runge}, existing solutions sequence
\(\{u_{1,j}\}_{j=1}^\infty\) and \(\{u_{2,j}^*\}_{j=1}^\infty\) such that
\begin{equation}\label{eq:knownA-runge-q}
  u_{1,j}|_\Omega=\psi_1+r_{1,j},
  \qquad
  u_{2,j}^*|_\Omega=\psi_2+r_{2,j}^*,
  \qquad
  r_{1,j}\to 0,\quad r_{2,j}^*\to 0
  \ \text{in }\widetilde{H}^s(\Omega).
\end{equation}
Inserting these into \eqref{eq:knownA-identity}, we obtain
\begin{equation}\label{eq:knownA-expansion-q}
\begin{aligned}
  0
  &=
  \int_\Omega (q_1-q_2)(\psi_1+r_{1,j})(\psi_2+r_{2,j}^*) \dd x \\
  &\quad+
  \int_\Omega (b_1-b_2)\cdot\nabla(\psi_1+r_{1,j})\,(\psi_2+r_{2,j}^*) \dd x.
\end{aligned}
\end{equation}
Since \(\psi_1=1\) on \(\operatorname{supp}(\psi_2)\), we have
\(\nabla\psi_1=0\) on \(\operatorname{supp}(\psi_2)\), and hence
\[
  \int_\Omega (b_1-b_2)\cdot\nabla\psi_1\,\psi_2 \dd x=0.
\]
Consequently, \eqref{eq:knownA-expansion-q} simplifies to
\begin{equation}\label{eq:knownA-main-q}
\begin{aligned}
  0
  &=
  \int_\Omega (q_1-q_2)\psi_1\psi_2 \dd x \\
  &\quad+
  \int_\Omega (q_1-q_2)\bigl(r_{1,j}\psi_2+\psi_1 r_{2,j}^*+r_{1,j}r_{2,j}^*\bigr) \dd x \\
  &\quad+
  \int_\Omega (b_1-b_2)\cdot\nabla r_{1,j}\,(\psi_2+r_{2,j}^*) \dd x.
\end{aligned}
\end{equation}
We estimate the error terms. By H\"older's inequality,
\begin{equation*}
\begin{aligned}
  &\Bigl|\int_\Omega (q_1-q_2)\bigl(r_{1,j}\psi_2+\psi_1 r_{2,j}^*+r_{1,j}r_{2,j}^*\bigr) \dd x\Bigr|\\
  &\le
  \|q_1-q_2\|_{L^\infty(\Omega)}
  \bigl(
    \|r_{1,j}\|_{L^2}\|\psi_2\|_{L^2}
    +
    \|\psi_1\|_{L^2}\|r_{2,j}^*\|_{L^2}
    +
    \|r_{1,j}\|_{L^2}\|r_{2,j}^*\|_{L^2}
  \bigr),
\end{aligned}
\end{equation*}
which tends to \(0\) as \(j\to\infty\). For the drift term, since \(b_1-b_2\in W^{1-s,\infty}(\Omega)\subset H^{1-s}(\Omega)\)
and \(\psi_2\in H^{1-s}(\Omega)\), the Kato--Ponce inequality gives
\[
  \|(b_1-b_2)\psi_2\|_{H^{1-s}(\Omega)}
  \le
  C\|b_1-b_2\|_{W^{1-s,\infty}(\Omega)}
  \|\psi_2\|_{H^{1-s}(\Omega)}.
\]
Using \(\nabla r_{1,j}\in H^{s-1}(\Omega)=(H^{1-s}(\Omega))^*\), we obtain
\[
  \Bigl|\int_\Omega (b_1-b_2)\cdot\nabla r_{1,j}\,\psi_2 \dd x\Bigr|
  \le
  \|b_1-b_2\|_{H^{s-1}(\Omega)}
  \|\nabla r_{1,j}\|_{H^{s-1}(\Omega)}
  \|\psi_2\|_{H^{1-s}(\Omega)}
  \to 0.
\]
as \(j\to\infty\), since \(\|\nabla r_{1,j}\|_{H^{s-1}}\le C\|r_{1,j}\|_{H^s}\to 0\).
The term involving \(r_{1,j}r_{2,j}^*\) is estimated similarly and also tends to \(0\).
Passing to the limit \(j\to\infty\) in \eqref{eq:knownA-main-q} therefore gives
\[
  \int_\Omega (q_1-q_2)\psi_1\psi_2 \dd x=0.
\]
Since \(\psi_1=1\) on \(\operatorname{supp}(\psi_2)\), this becomes
\[
  \int_\Omega (q_1-q_2)\psi_2 \dd x=0.
\]
As \(\psi_2\in C_c^\infty(\Omega)\) is arbitrary, we can obtain \(q_1=q_2\) in \(\Omega\).

\emph{Step 2: Recovery of \(b\).}
With \(q_1=q_2\), identity \eqref{eq:knownA-identity} reduces to
\begin{equation}\label{eq:knownA-drift}
  \int_\Omega (b_1-b_2)\cdot\nabla u_1\, u_2^* \dd x=0.
\end{equation}
Fix \(k\in\{1,\dots,n\}\) and let \(\psi_2\in C_c^\infty(\Omega)\) be arbitrary.
Choose a cutoff function \(\psi_{x_k}\in C_c^\infty(\Omega)\) such that
\(\psi_{x_k}(x)=x_k\) on \(\operatorname{supp}(\psi_2)\).
By the Runge approximation, there exist sequences of solutions
\(\{u_{1,j}\}_{j=1}^\infty\) and \(\{u_{2,j}^*\}_{j=1}^\infty\) satisfying
\begin{equation}\label{eq:knownA-runge-b}
  u_{1,j}|_\Omega=\psi_{x_k}+r_{1,j},
  \qquad
  u_{2,j}^*|_\Omega=\psi_2+r_{2,j}^*,
  \qquad
  r_{1,j}\to 0,\quad r_{2,j}^*\to 0
  \ \text{in }\widetilde{H}^s(\Omega).
\end{equation}
Substituting into \eqref{eq:knownA-drift} yields
\begin{equation}\label{eq:knownA-expansion-b}
\begin{aligned}
  0
  &=
  \int_\Omega (b_1-b_2)\cdot\nabla\psi_{x_k}\,\psi_2 \dd x \\
  &\quad+
  \int_\Omega (b_1-b_2)\cdot\nabla\psi_{x_k}\,r_{2,j}^* \dd x \\
  &\quad+
  \int_\Omega (b_1-b_2)\cdot\nabla r_{1,j}\,(\psi_2+r_{2,j}^*) \dd x.
\end{aligned}
\end{equation}
Since \(\nabla\psi_{x_k}=e_k\) on \(\operatorname{supp}(\psi_2)\), the first term equals
\[
 \int_\Omega (b_1-b_2)\cdot\nabla\psi_{x_k}\,\psi_2 \dd x= \int_\Omega (b_1-b_2)_k\,\psi_2 \dd x.
\]
The error terms tend to zero as \(j\to\infty\), by the same argument as in Step 1.
Passing to the limit in \eqref{eq:knownA-expansion-b} gives
\[
  \int_\Omega (b_1-b_2)_k\,\psi_2 \dd x=0.
\]
Since \(k\in\{1,\dots,n\}\) and \(\psi_2\in C_c^\infty(\Omega)\) are arbitrary, we conclude
\((b_1-b_2)_k=0\) for every \(k\), hence \(b_1=b_2\) in \(\Omega\).
This completes the proof.
\end{proof}
\subsection{Case II: Unknown magnetic potential}\label{sec:unknownA}
\begin{proof}[Proof of Theorem \ref{thm:unknownA}]
Without loss of generality, we may assume \(W_1\cap W_2=\emptyset\); indeed,
if \(W_1\cap W_2\neq\emptyset\), we replace \(W_i\) by smaller non-empty open
subsets \(W_i'\subset W_i\) with \(W_1'\cap W_2'=\emptyset\), and the hypothesis
remains valid on \(W_1',W_2'\).

Let \(f_i\in C_c^\infty(W_i)\) and let \(u_1\in H^s(\Rb^n)\) solve
\[
  (-\Delta)_{A_1}^s u_1+b_1\cdot\nabla u_1+q_1u_1=0\quad\text{in }\Omega,
  \qquad u_1-f_1\in\widetilde{H}^s(\Omega),
\]
while \(u_2^*\in H^s(\Rb^n)\) solves the adjoint equation
\[
  (-\Delta)_{A_2}^s u_2^*-\nabla\cdot(b_2u_2^*)+q_2u_2^*=0\quad\text{in }\Omega,
  \qquad u_2^*-f_2\in\widetilde{H}^s(\Omega).
\]
Since the DN maps computed on \(f\in C_c^\infty(W_1)\) are equal after restriction on \(W_2\), the Alessandrini identity \eqref{eq:alessandrini} yields
\begin{equation}\label{eq:unknownA-Aless}
\begin{aligned}
0
&=2\left\langle\int_{\Rb^n}\bigl((A_1)_{a\parallel}-(A_2)_{a\parallel}\bigr)
\cdot\nabla^s u_1\dd y,\;u_2^*\right\rangle\\
&\quad+\langle(Q_1-Q_2)u_1,u_2^*\rangle
+\int_\Omega(b_1-b_2)\cdot\nabla u_1\,u_2^*\dd x.
\end{aligned}
\end{equation}

Since \(u_i\) is supported in \(\Omega\cup W_i\) for \(i=1,2\) and
\((\Omega\cup W_1)\cap(\Omega\cup W_2)=\Omega\) (as \(W_1\cap W_2=\emptyset\)),
and since \((A_i,q_i)\in\mathcal P\) so that \(\supp(A_i)\subseteq\Omega^2\),
we have
\begin{equation}\label{eq:unknownA-23}
\begin{aligned}
0
&=2\int_\Omega(u_2^*|_\Omega)(x)\int_\Omega(\sigma_1-\sigma_2)(x,y)|\alpha|^2
\bigl((u_1|_\Omega)(x)-(u_1|_\Omega)(y)\bigr)\dd y\dd x\\
&\quad+\int_\Omega u_1|_\Omega\,u_2^*|_\Omega\,(Q_1-Q_2)\dd x
+\int_\Omega(b_1-b_2)\cdot\nabla(u_1|_\Omega)\,u_2^*|_\Omega\dd x,
\end{aligned}
\end{equation}
where we used Lemma \ref{lem:extension2} and \(|\alpha|^2=C_{n,s}|x-y|^{-n-2s}\).

Fix \(g\in C_c^\infty(\Omega)\) and \(f\in C_c^\infty(\Omega)\). By Lemma~\ref{lem:runge}(b)
and its adjoint version, there exist sequences of solutions
\(u_1^{(k)}=P_{A_1,b_1,q_1}f_1^{(k)}\) and
\(u_2^{(k,*)}=P_{A_2,b_2,q_2}^*f_2^{(k)}\) with \(f_i^{(k)}\in C_c^\infty(W_i)\) such that
\[
  u_1^{(k)}|_\Omega\to g,
  \qquad
  u_2^{(k,*)}|_\Omega\to f
  \qquad\text{in }\widetilde{H}^s(\Omega)
\]
as \(k\to\infty\). Inserting these into \eqref{eq:unknownA-23} and passing to
the limit \(k\to\infty\), we obtain
\begin{equation}\label{eq:unknownA-general}
\begin{aligned}
0
&=2\int_\Omega(\sigma_1-\sigma_2)(x,y)|\alpha|^2
\bigl(g(x)-g(y)\bigr)\dd y
+g(x)(Q_1-Q_2)(x)\\
&\quad+(b_1-b_2)(x)\cdot\nabla g(x)
\qquad\text{for a.e. }x\in\Omega.
\end{aligned}
\end{equation}
Since \(f,g\in C_c^\infty(\Omega)\) are arbitrary and \(C_c^\infty(\Omega)\) is
dense in \(L^2(\Omega)\), identity \eqref{eq:unknownA-general} holds for all
\(g\in C_c^\infty(\Omega)\) and a.e. \(x\in\Omega\).

\emph{Step 1: Recovery of \((A_1)_{a\parallel}\).}
Fix \(x_0\in\Omega\) and \(\psi\in C_c^\infty(\Omega)\), and define
\[
  g(y):=\psi(y)\,e^{-1/|x_0-y|}\,|x_0-y|^{n+2s}.
\]
Then \(g\in C_c^\infty(\Omega)\), \(g(x_0)=0\) and \(\nabla g(x_0)=0\).
Substituting this \(g\) into \eqref{eq:unknownA-general} and evaluating at
\(x=x_0\) yields
\[
  -\int_\Omega(\sigma_1-\sigma_2)(x_0,y)|\alpha(x_0,y)|^2g(y)\dd y=0.
\]
Since \(|\alpha(x_0,y)|^2=C_{n,s}|x_0-y|^{-n-2s}\), this simplifies to
\[
  \int_\Omega(\sigma_1-\sigma_2)(x_0,y)\,\psi(y)\,e^{-1/|x_0-y|}\dd y=0.
\]
Since \(\psi\in C_c^\infty(\Omega)\) is arbitrary and the exponential factor is
strictly positive for \(y\neq x_0\), we deduce
\(\sigma_1(x_0,y)=\sigma_2(x_0,y)\) for a.e. \(y\in\Omega\).
As \(x_0\in\Omega\) is arbitrary and \(\sigma_1,\sigma_2\) are symmetric with
\(\sigma_1=\sigma_2=1\) on \(\Rb^{2n}\setminus\Omega^2\), this implies
\(\sigma_1=\sigma_2\) on \(\Rb^{2n}\). Therefore, we get
\[
  (A_1)_{a\parallel}=(A_2)_{a\parallel}.
\]

\emph{Step 2: Recovery of \(Q\).}
With \(\sigma_1=\sigma_2\), identity \eqref{eq:unknownA-general} reduces to
\[
  g(x)(Q_1-Q_2)(x)+(b_1-b_2)(x)\cdot\nabla g(x)=0
  \qquad\text{for a.e. }x\in\Omega,\ \forall g\in C_c^\infty(\Omega).
\]
Fix \(x\in\Omega\) and choose \(g\in C_c^\infty(\Omega)\) equal to \(1\) in a
neighborhood of \(x\); then \(\nabla g(x)=0\) and we obtain \((Q_1-Q_2)(x)=0\).
As \(x\in\Omega\) is arbitrary, \(Q_1=Q_2\) in \(\Omega\).

\emph{Step 3: Recovery of \(b\).}
With \(Q_1=Q_2\), the identity further reduces to
\[
  (b_1-b_2)(x)\cdot\nabla g(x)=0
  \qquad\text{for a.e. }x\in\Omega,\ \forall g\in C_c^\infty(\Omega).
\]
Fix \(k\in\{1,\dots,n\}\), choose a cutoff function \(\psi\in C_c^\infty(\Omega)\)
with \(\psi=1\) on \(\operatorname{supp}(b_1-b_2)\), and take \(g(x)=\psi(x)x_k\).
Then \(\nabla g=e_k\) on \(\operatorname{supp}(b_1-b_2)\), and we obtain
\((b_1-b_2)_k=0\). Since \(k\) is arbitrary, \(b_1=b_2\) in \(\Omega\).

Combining \((A_1)_{a\parallel}=(A_2)_{a\parallel}\), \(b_1=b_2\) and \(Q_1=Q_2\),
we can obtain \((A_1,b_1,q_1)\sim(A_2,b_2,q_2)\).
\end{proof}

\begin{remark}
Theorem~\ref{thm:unknownA} shows that \((A,b,q)\) cannot be uniquely determined
when \(A\) is unknown. The gauge freedom provides infinitely many equivalent
triples producing the same DN map. However, the drift \(b\) is a gauge invariant
and is uniquely recovered.
\end{remark}

\section{Higher-order Runge approximation}\label{sec:higher-order}
In this section, we first review the properties of high-order Runge approximation and provide the proof of Theorem \ref{thm:runge}(b). The proof structure of high-order Runge approximation is similar to the one discussed in Section \ref{Runge}. The reasoning of this argument is basically consistent with the proof in Section 7 of \cite{GSU-2016}, and for the sake of self-consistency, we will restate its content below, which is applicable to the fractional magnetic Schr\"odinger equation with a drift term.

Throughout this section we fix \(s\in(1/2,1)\) whose boundary is \(C^\infty\)-smooth. Concerning the coefficients,
we require
\[
  A\in C_c^\infty(\Omega^2),
  \qquad
  b\in C_c^\infty(\Omega)^n,
  \qquad
  q\in C_c^\infty(\Omega),
  \qquad
  \supp(b),\supp(q)\Subset\Omega.
\]

We consider the following function space
\[
  \mathcal{E}^s(\overline{\Omega})
  :=
  e_\Omega\bigl(d(x)^s C^\infty(\overline{\Omega})\bigr),
\]
where \(e_\Omega\) stands for extension by zero from \(\Omega\) to \(\Rb^n\),
while \(d(x)\) denotes a smooth positive function on \(\overline{\Omega}\) that
coincides with \(\operatorname{dist}(x,\partial\Omega)\) in a neighbourhood of
\(\partial\Omega\).

Given \(\mu>s-\frac12\), we further employ the scale of Banach spaces
\(H^{s(\mu)}(\overline{\Omega})\) introduced by Grubb \cite{Grubb2015}.
These spaces are designed precisely for functions \(u\) satisfying
\(u=0\) on \(\Omega_e\) and such that
\[
  r_\Omega\bigl((-\Delta)_A^s+b\cdot\nabla+q\bigr)u
  \in H^{\mu-2s}(\Omega),
\]
where \(r_\Omega\) is the restriction operator \(r_\Omega u:=u|_\Omega\).
The key feature of \(H^{s(\mu)}(\overline{\Omega})\) is that, on this scale,
the magnetic operator with drift acts as an isomorphism onto
\(H^{\mu-2s}(\Omega)\); this will be made precise in
Lemma~\ref{prop:magnetic-mapping} below.

To derive the required higher-order approximation properties, we recall the following result.

\begin{proposition}\cite{Grubb2015}\label{prop:Grubb-spaces}
Let \(\mu>s-\frac12\) and let \(\Omega\subset\Rb^n\) be a smooth bounded domain.
There exists a Banach space \(H^{s(\mu)}(\overline{\Omega})\) with the following
properties:
\begin{enumerate}
  \item[(a)] \(H^{s(\mu)}(\overline{\Omega})\subset H^{s-1/2}_{\overline{\Omega}}\)
  with continuous inclusion;
  \item[(b)] \(H^{s(\mu)}(\overline{\Omega})=H^{\mu}_{\overline{\Omega}}\) for
  \(\mu\in(s-\frac12,s+\frac12)\);
  \item[(c)] \(r_\Omega((-\Delta)^s)\) is a homeomorphism from
  \(H^{s(\mu)}(\overline{\Omega})\) onto \(H^{\mu-2s}(\Omega)\);
  \item[(d)] \(H^{\mu}_{\overline{\Omega}}\subset H^{s(\mu)}(\overline{\Omega})
  \subset H^{\mu}_{\mathrm{loc}}(\Omega)\) with continuous inclusions, and
  multiplication by any \(\chi\in C_c^\infty(\Omega)\) is bounded from
  \(H^{s(\mu)}(\overline{\Omega})\) to \(H^\mu(\Omega)\);
  \item[(e)] the space \(\mathcal{E}^s(\overline{\Omega}):=\bigcap_{\mu>s-1/2}H^{s(\mu)}(\overline{\Omega})\)
  is dense in every \(H^{s(\mu)}(\overline{\Omega})\).
\end{enumerate}
\end{proposition}

A key technical point in the higher-order approximation argument is the
behaviour of \((-\Delta)_A^s+b\cdot\nabla+q\) on test functions supported
in \(\Omega_e\). In contrast to the non-magnetic case, where the drift and
potential terms vanish trivially on \(\Omega\), the magnetic contribution
requires a separate verification. The following lemma provides this
reduction and will be invoked in the duality argument of
Lemma~\ref{lem:higher-order-runge}.

\begin{lemma}\label{lem:A-vanishes-Omega}
Let \(A\in C_c^\infty(\Omega^2)\), \(b,q\in C_c^\infty(\Omega)\), and let
\(f\in C_c^\infty(\Omega_e)\). Then for every \(x\in\Omega\),
\[
  \int_{\Rb^n}A_{a\parallel}(x,y)\cdot\nabla^s f(x,y)\dd y=0,
  \qquad
  b(x)\cdot\nabla f(x)+q(x)f(x)=0.
\]
In particular,
\[
  r_\Omega\bigl((-\Delta)_A^s f+b\cdot\nabla f+qf\bigr)
  =
  r_\Omega((-\Delta)^s f).
\]
\end{lemma}

\begin{proof}
Since \(\supp(A)\subseteq\Omega^2\), we have \(A(x,y)=0\) whenever \(x\notin\Omega\)
or \(y\notin\Omega\). Fix \(x\in\Omega\). The integrand in
\(\int_{\Rb^n}A_{a\parallel}(x,y)\cdot\nabla^s f(x,y)\dd y\) is supported in
\(y\in\Omega\), where \(f(y)=0\) (as \(f\in C_c^\infty(\Omega_e)\)). Hence
\[
  \nabla^s f(x,y)=(f(x)-f(y))\alpha(x,y)
\]
vanishes identically on the support of \(A_{a\parallel}(x,\cdot)\), because
\(f(x)=0\) (as \(x\in\Omega\)) and \(f(y)=0\) (as \(y\in\Omega\)). Therefore
the integral vanishes. Similarly, \(b,q\in C_c^\infty(\Omega)\) and \(f=0\) on
\(\Omega\) give \(b(x)\cdot\nabla f(x)+q(x)f(x)=0\).
\end{proof}

We now prove that the magnetic operator with drift maintains the mapping properties of the fractional Laplacian at the Grubb scale, that is, Proposition \ref{prop:Grubb-spaces}(c) still holds.

\begin{lemma}\label{prop:magnetic-mapping}
Under the assumptions of Theorem~\ref{thm:runge}(b), for every
\(\mu>s-\frac12\) the operator
\[
  \mathcal{T}_{A,b,q}:=r_\Omega\bigl((-\Delta)_A^s+b\cdot\nabla+q\bigr)
\]
is a homeomorphism from \(H^{s(\mu)}(\overline{\Omega})\) onto
\(H^{\mu-2s}(\Omega)\).
\end{lemma}

\begin{proof}
By the decomposition \eqref{eq22}, we write
\[
  \mathcal{T}_{A,b,q}=\mathcal{T}_0+\mathcal{F}_A+\mathcal{F}_{b,q},
\]
where
\[
  \mathcal{T}_0=r_\Omega((-\Delta)^s),\qquad
  \mathcal{F}_A u=r_\Omega\Bigl(2\int_{\Rb^n}A_{a\parallel}\cdot\nabla^s u\dd y+Qu\Bigr),
  \qquad
  \mathcal{F}_{b,q}u=r_\Omega(b\cdot\nabla u+qu),
\]
with \(Q=q+\int_{\Rb^n}|A|^2\dd y+(\nabla\cdot)^s A_{s\parallel}\).
By Proposition~\ref{prop:Grubb-spaces}(c), \(\mathcal{T}_0\) is a homeomorphism
from \(H^{s(\mu)}(\overline{\Omega})\) onto \(H^{\mu-2s}(\Omega)\).

\emph{Boundedness of \(\mathcal{F}_{b,q}\).} For \(u\in H^{s(\mu)}(\overline{\Omega})\),
Proposition~\ref{prop:Grubb-spaces}(d) combined with \(b,q\in C_c^\infty(\Omega)\) gives
\[
  \|\mathcal{F}_{b,q}u\|_{H^{\mu-1}(\Omega)}
  \le C(\|b\|_{C^\infty}+\|q\|_{C^\infty})\|u\|_{H^{s(\mu)}(\overline{\Omega})}.
\]
Since \(2s>1\) and \(\mu-1>\mu-2s\), we have \(H^{\mu-1}(\Omega)\subset
H^{\mu-2s}(\Omega)\) with continuous inclusion, and the embedding is actually
compact because \(\mu-1>\mu-2s\). Therefore, when \(\mathcal{F}_{b,q}\) is an operator on \(H^{\mu-2s}(\Omega)\), it is bounded and compact from \(H^{s(\mu)}(\overline{\Omega})\) to \(H^{\mu-2s}(\Omega)\).

\emph{Boundedness of \(\mathcal{F}_A\).} Since \(A\in C_c^\infty(\Omega^2)\), the function
\(Q_0:=\int_{\Rb^n}|A|^2\dd y\in C_c^\infty(\Omega)\) and
\((\nabla\cdot)^s A_{s\parallel}\in C_c^\infty(\Omega)\); hence
\(Q\in C_c^\infty(\Omega)\). The multiplication operator \(u\mapsto Qu\) is
therefore bounded from \(H^{s(\mu)}(\overline{\Omega})\) to \(H^\mu(\Omega)\subset
H^{\mu-s}(\Omega)\). For the integral term, set
\[
  K_A u(x):=2\int_{\Rb^n}A_{a\parallel}(x,y)\cdot\alpha(x,y)\,
  \bigl(u(x)-u(y)\bigr)\dd y.
\]
Since \(A\in C_c^\infty(\Omega^2)\) and \(\alpha\) is homogeneous of degree
\(-n/2-s+1\), the kernel \(k_A(x,y):=A_{a\parallel}(x,y)\cdot\alpha(x,y)\)
is smooth off the diagonal, supported in \(\Omega^2\), and satisfies the
symbol-type estimates
\[
  |\partial_x^\beta\partial_y^\gamma k_A(x,y)|
  \le C_{\beta,\gamma}|x-y|^{-n/2-s+1-|\gamma|}
  \qquad \text{for }|x-y|\le 1.
\]
For \(s\in(1/2,1)\) and \(n\ge1\), one has \(-n/2-s+1>-n\), hence the kernel
\(k_A\) is locally integrable. Splitting
\[
  K_Au(x)
  =
  2u(x)\int_{\Rb^n}k_A(x,y)\dd y
  -
  2\int_{\Rb^n}k_A(x,y)u(y)\dd y
  =:
  M_A(x)u(x)+L_Au(x),
\]
the first term is a multiplication operator with smooth coefficient
\(M_A\in C_c^\infty(\Omega)\), hence bounded from \(H^{s(\mu)}(\overline{\Omega})\)
to \(H^\mu(\Omega)\).
For \(L_A\), we estimate the kernel
\[
  k_A(x,y):=A_{a\parallel}(x,y)\cdot\alpha(x,y).
\]
 Since \(A\in C_c^\infty(\Omega^2)\) and \(A_{a\parallel}\) is
antisymmetric, we have \(A_{a\parallel}(x,x)=0\) and hence
\[
  |A_{a\parallel}(x,y)|\le C|x-y|.
\]
Together with \(|\alpha(x,y)|=C_{n,s}^{1/2}|x-y|^{-n/2-s}\), this yields the
kernel bound
\begin{equation}\label{eq:kernel-bound}
  |k_A(x,y)|\le C|x-y|^{-n/2-s+1}.
\end{equation}
The kernel \(k_A\) is smooth off the diagonal and supported in \(\Omega^2\).
By the theory of fractional integral operators (equivalently, by the
pseudodifferential calculus of \cite[Chapter 5]{Grubb2015}), the operator
\(L_A\) with kernel \eqref{eq:kernel-bound} is smoothing of order
\[
  s_0:=n/2+1-s>0,
\]
i.e.\ it maps \(H^p(\Rb^n)\) boundedly into \(H^{p+s_0}(\Rb^n)\) for every
\(p\in\Rb\). Since \(s_0>0\) and \(H^{p+s_0}\subset H^{p-s}\) for every
\(p\in\Rb\), we obtain in particular
\[
  \|L_Au\|_{H^{\mu-s}(\Omega)}
  \le C\|u\|_{H^\mu(\Rb^n)}
  \qquad\text{for all }u\in H^\mu(\Rb^n).
\]
Using \(H^{s(\mu)}(\overline{\Omega})\subset H^\mu_{\mathrm{loc}}(\Omega)\)
from Proposition~\ref{prop:Grubb-spaces}(d), we conclude
\[
  \|L_Au\|_{H^{\mu-s}(\Omega)}
  \le C\|u\|_{H^{s(\mu)}(\overline{\Omega})}
  \qquad\text{for all }u\in H^{s(\mu)}(\overline{\Omega}).
\]
Since \(\mu-s>\mu-2s\), the embedding \(H^{\mu-s}(\Omega)\hookrightarrow
H^{\mu-2s}(\Omega)\) is compact. Combining the estimates for \(M_A\) and
\(L_A\) gives the boundedness of \(\mathcal{F}_A\), and its compactness as an
operator on \(H^{\mu-2s}(\Omega)\).

By the above, \(\mathcal{F}_A+\mathcal{F}_{b,q}\) is
bounded from \(H^{s(\mu)}(\overline{\Omega})\) to \(H^{\mu-2s}(\Omega)\) and
compact when regarded as an operator \(H^{\mu-2s}(\Omega)\to H^{\mu-2s}(\Omega)\).
Therefore \(\mathcal{T}_{A,b,q}=\mathcal{T}_0(I+\mathcal{T}_0^{-1}
(\mathcal{F}_A+\mathcal{F}_{b,q}))\) is Fredholm of index zero from
\(H^{s(\mu)}(\overline{\Omega})\) to \(H^{\mu-2s}(\Omega)\).

Suppose \(\mathcal{T}_{A,b,q}u=0\) for some
\(u\in H^{s(\mu)}(\overline{\Omega})\subset H^{s-1/2}_{\overline{\Omega}}\).
Then \(u\in\widetilde{H}^s(\Omega)\) and
\(((-\Delta)_A^s+b\cdot\nabla+q)u=0\) in \(\Omega\). By the eigenvalue
condition, \(u\equiv0\). Hence \(\mathcal{T}_{A,b,q}\) is injective, and being
Fredholm of index zero, it is an isomorphism.
\end{proof}

Based on these properties and following the argument in \cite{GSU-2016}, we have proved the high-order approximation property of the solution of equation \eqref{eq11} in $\mathcal{E}^s(\overline{\Omega})$.
\begin{lemma}\label{lem:higher-order-runge}
If the conditions in Theorem \ref{thm:runge}(b) are satisfied, let \(P_{A,b,q}\) denote
the Poisson operator from Definition~\ref{def:DN} and set
\[
  \mathcal{D}:=\bigl\{e_\Omega(r_\Omega P_{A,b,q}f):f\in C_c^\infty(W)\bigr\}.
\]
Then \(\mathcal{E}^s(\overline{\Omega})\) forms a Fréchet space equipped with the topology generated by \(\{\|\cdot\|_{H^{s(k)}}\}_{k\in\mathbb{N}}\), and \(\mathcal{D}\) is dense in this space.
\end{lemma}

\begin{proof}
Fix \(f\in C_c^\infty(W)\subset C_c^\infty(\Omega_e)\) and set
\(u_f:=P_{A,b,q}f\). Then \(u_f-f\in\widetilde{H}^s(\Omega)\) and
\(\mathcal{T}_{A,b,q}(u_f-f)=-\mathcal{T}_{A,b,q}f\). By
Lemma~\ref{lem:A-vanishes-Omega},
\[
  \mathcal{T}_{A,b,q}f=r_\Omega((-\Delta)^s f)\in C^\infty(\overline{\Omega}),
\]
and \((-\Delta)^s f\in H^{\mu-2s}(\Omega)\) for every \(\mu>s-\frac12\). By
Lemma~\ref{prop:magnetic-mapping}, \(u_f-f\in H^{s(\mu)}(\overline{\Omega})\)
for every \(\mu>s-\frac12\), hence \(u_f-f\in\mathcal{E}^s(\overline{\Omega})\).
Consequently \(e_\Omega(r_\Omega u_f)=e_\Omega(r_\Omega(u_f-f))\in\mathcal{E}^s(\overline{\Omega})\).

Let \(\mathcal{L}\) be a continuous
linear functional on \(\mathcal{E}^s(\overline{\Omega})\) satisfying
\[
  \mathcal{L}\bigl(e_\Omega(r_\Omega P_{A,b,q}f)\bigr)=0
  \qquad\text{for all }f\in C_c^\infty(W).
\]
By the Hahn-Banach theorem, it suffices to show that \(\mathcal{L}=0\).

Since \(\mathcal{L}\) is continuous
on the Fr\'echet space \(\mathcal{E}^s(\overline{\Omega})\), there exist
\(\ell\in\Nb\) and \(C>0\) such that
\[
  |\mathcal{L}(u)|\le C\|u\|_{H^{s(\ell)}(\overline{\Omega})}
  \qquad\text{for all }u\in\mathcal{E}^s(\overline{\Omega}).
\]
By Proposition~\ref{prop:Grubb-spaces}(e), \(\mathcal{E}^s(\overline{\Omega})\)
is dense in \(H^{s(\ell)}(\overline{\Omega})\); hence \(\mathcal{L}\) extends
uniquely to a bounded linear functional
\(\widetilde{\mathcal{L}}\in(H^{s(\ell)}(\overline{\Omega}))^*\).

Consider the isomorphism
\(\mathcal{T}_{A,b,q}:H^{s(\ell)}(\overline{\Omega})\to H^{\ell-2s}(\Omega)\)
from Lemma~\ref{prop:magnetic-mapping}. Its adjoint
\[
  \mathcal{T}_{A,b,q}^*:(H^{\ell-2s}(\Omega))^*\to(H^{s(\ell)}(\overline{\Omega}))^*
\]
is again an isomorphism. By \cite{CekicLinRuland2020},
\[
  (H^{\ell-2s}(\Omega))^*
  =
  \widetilde{H}^{-\ell+2s}(\Omega)
  =
  H^{-\ell+2s}_{\overline{\Omega}},
\]
so we may identify \(\mathcal{T}_{A,b,q}^*\) with the map
\(v\mapsto(w\mapsto(v,\mathcal{T}_{A,b,q}w))\) for
\(v\in H^{-\ell+2s}_{\overline{\Omega}}\). Choose
\(v\in H^{-\ell+2s}_{\overline{\Omega}}\) with
\(\mathcal{T}_{A,b,q}^*v=\widetilde{\mathcal{L}}\).

For \(f\in C_c^\infty(W)\),
since \(P_{A,b,q}f\) solves the homogeneous equation in \(\Omega\), we have
\(\mathcal{T}_{A,b,q}(P_{A,b,q}f)=0\), and since \(P_{A,b,q}f-f\in\widetilde{H}^s(\Omega)\),
we obtain
\begin{equation}\label{eq:higher-order-duality}
\begin{aligned}
  0
  &=\mathcal{L}\bigl(e_\Omega(r_\Omega P_{A,b,q}f)\bigr)
  =\widetilde{\mathcal{L}}(P_{A,b,q}f-f)
  =\mathcal{T}_{A,b,q}^*v(P_{A,b,q}f-f)\\
  &=(v,\mathcal{T}_{A,b,q}(P_{A,b,q}f-f))
  =-(v,\mathcal{T}_{A,b,q}f).
\end{aligned}
\end{equation}
By Lemma~\ref{lem:A-vanishes-Omega}, \(\mathcal{T}_{A,b,q}f=r_\Omega((-\Delta)^s f)\),
so
\[
  (v,(-\Delta)^s f)=0\qquad\text{for all }f\in C_c^\infty(W).
\]
By the self-adjointness of \((-\Delta)^s\), this can be rewritten as
\[
  \bigl((-\Delta)^s v,f\bigr)=0\qquad\text{for all }f\in C_c^\infty(W).
\]
Hence \((-\Delta)^s v=0\) in \(W\) in the sense of distributions.

Since
\(v\in H^{-\ell+2s}_{\overline{\Omega}}\subset H^{-\ell+2s}(\Rb^n)\) and
\(v=0\) on \(\Omega_e\supset W\), we have
\[
  v=(-\Delta)^s v=0\quad\text{in }W.
\]
Consequently \(\widetilde{\mathcal{L}}=0\) and hence \(\mathcal{L}=0\), which
proves the density of \(\mathcal{D}\).
\end{proof}

We next verify the Runge approximation property in the setting of higher regularity. T
\begin{proof}[Proof of Theorem~\ref{thm:runge}(b)]
Let \(g\in C^\infty(\overline{\Omega})\) and define
\(h:=e_\Omega d(x)^s g\in\mathcal{E}^s(\overline{\Omega})\). By
Lemma~\ref{lem:higher-order-runge}, there exists a sequence
\(\{f_j\}_{j\in\Nb}\subset C_c^\infty(W)\) such that
\[
  e_\Omega(r_\Omega P_{A,b,q}f_j)\to h
  \qquad\text{in }\mathcal{E}^s(\overline{\Omega})
\]
as \(j\to\infty\).
Next, consider the map
\[
  \mathcal{M}:C^\infty(\overline{\Omega})\to\mathcal{E}^s(\overline{\Omega}),
  \qquad
  \mathcal{M}g:=e_\Omega d(x)^s g.
\]
The higher-order approximation will follow once we know that \(\mathcal{M}\) is
a homeomorphism, since its inverse is then necessarily given by
\(\mathcal{M}^{-1}=d(x)^{-s}r_\Omega\).

Now, Given that \(\mathcal{M}\) is a bijective linear map between Fr\'echet spaces. Its
graph is closed, if \(g_j\to g\) in \(C^\infty(\overline{\Omega})\) and
\(\mathcal{M}g_j\to h\) in \(\mathcal{E}^s(\overline{\Omega})\), then
\(g_j\to g\) in \(L^\infty(\overline{\Omega})\) and
\(\mathcal{M}g_j\to\mathcal{M}g\) in \(L^2(\Rb^n)\), while also
\(\mathcal{M}g_j\to h\) in \(L^2(\Rb^n)\); by uniqueness of the limit,
\(h=\mathcal{M}g\) in the sense of distributions. As a consequence, the closed graph theorem together with the open mapping theorem implies that \(\mathcal{M}\) is a homeomorphism. Applying
\(\mathcal{M}^{-1}=d(x)^{-s}r_\Omega\) to the convergence above, we obtain
\[
  d(x)^{-s}r_\Omega P_{A,b,q}f_j\to g
  \qquad\text{in }C^\infty(\overline{\Omega}),
\]
and since convergence in \(C^\infty(\overline{\Omega})\) implies convergence in
\(C^k(\Omega)\) for every \(k\in\Nb\), this yields the desired conclusion.
\end{proof}

\section{Numerical method}\label{sec:numerical}

In this section, we establish a numerical framework for inverting the drift term \(b\) and the scalar potential \(q\) from external Dirichlet-Neumann quantitative measurement data, assuming that the magnetic potential \(A\) is known. By discretizing the forward operator using the finite difference method on the truncated domain, the unknown parameters \(b\) and \(q\) are parameterized as low-dimensional bases, and a radial basis function (RBF) surrogate model is constructed. The inversion process employs a smooth Tikhonov-TV regularization method, thereby avoiding the repeated solution of partial differential equations during the optimization process.

\subsection{Truncation and finite-difference discretization}\label{sec:num-fd}

Let \(\Omega=(-L,L)\) and let \(R>L\) be a truncation parameter. For an
exterior datum \(g\in C_c^\infty(W)\) with \(W\subset\Omega_e\), we consider
\begin{equation}\label{eq:num-truncated}
\begin{cases}
\bigl((-\Delta)_A^s+b\,\partial_x+q\bigr)u_R=0, & x\in(-L,L),\\[2pt]
u_R=f, & x\in(-R,R)\setminus(-L,L),\\[2pt]
u_R=0, & |x|\ge R,
\end{cases}
\end{equation}
where \(A\) is fixed. Let \(N\in\Nb\), \(h=2L/N\), and let
\[
x_i=-L+ih,\qquad i=1,\dots,N-1,
\]
be the interior grid points. We write
\[
\mathbf U_R=(u_R(x_1),\dots,u_R(x_{N-1}))^\top\in\Rb^{N-1}.
\]
The exterior grid points in \((-R,R)\setminus(-L,L)\) carry the datum \(f\),
and the points outside \((-R,R)\) are set to zero.

The drift and the potential are parameterized by
\begin{equation}\label{eq:num-bq-param}
b(x)=\sum_{k=1}^{N_b}\beta_k\psi_k(x),
\qquad
q(x)=\sum_{k=1}^{N_q}\theta_k\chi_k(x),
\end{equation}
where \(\psi_k,\chi_k\in C_c^\infty(\Omega)\). The parameter vector is
\[
\boldsymbol\xi=(\beta,\theta)\in\Rb^M,
\qquad
M=N_b+N_q.
\]

Writing
\[
(-\Delta)^s u_R(x_i)
=-c_{1,s}\int_0^\infty
\frac{u_R(x_i-\xi)-2u_R(x_i)+u_R(x_i+\xi)}{\xi^{1+2s}}\dd\xi
\]
and introducing  \(\rho\in(2s,2]\), the integral is approximated by a
weighted sum over the grid, yielding a symmetric Toeplitz-like matrix
\(\mathbf L^s\in\Rb^{(N-1)\times(N-1)}\) acting on \(\mathbf U_R\)\cite{Li2025}.

Since \(A\) is known and fixed,
we assume that its admissible part is given in the form
\[
A_{a\parallel}(x,y)=\alpha(x,y)\bigl(\gamma(x)+\gamma(y)\bigr),
\]
for some known function \(\gamma\in C_c^\infty(\Omega)\), where \(\alpha(x,y)\) denotes the kernel of the fractional gradient,
\(\nabla^s u(x,y)=(u(x)-u(y))\alpha(x,y)\), whose modulus satisfies \(|\alpha(x,y)|^2=\frac{c_{1,s}}2|x-y|^{-1-2s}\) in one dimension. By \eqref{eq22},
\[
\bigl((-\Delta)_A^s-(-\Delta)^s-Q\bigr)u(x)
=2\int_{\Rb}A_{a\parallel}(x,y)\cdot\nabla^s u(x,y)\dd y,
\]
and, since \(A_{a\parallel}\) is parallel to \(\alpha\), the product becomes
\[
2A_{a\parallel}\cdot\nabla^s u
=2\,|\alpha(x,y)|^2\bigl(\gamma(x)+\gamma(y)\bigr)\bigl(u(x)-u(y)\bigr).
\]
So
\[
2\int_{\Rb}A_{a\parallel}\cdot\nabla^s u\dd y
=
c_{1,s}\int_{\Rb}|x-y|^{-1-2s}
\bigl(\gamma(x)+\gamma(y)\bigr)\bigl(u(x)-u(y)\bigr)\dd y.
\]
This yields a \emph{fixed} matrix
\(\mathbf M_A\in\Rb^{(N-1)\times(N-1)}\) with off-diagonal
entries
\begin{equation}\label{eq:num-MA}
(\mathbf M_A)_{ij}
=
-\kappa_{ij}\,(\gamma_i+\gamma_j),
\qquad i\neq j,
\end{equation}
where
\[
\kappa_{ij}:=w_{ij}\,c_{1,s}\,|x_i-x_j|^{-1-2s},
\qquad
\gamma_i=\gamma(x_i),
\]
and \(w_{ij}\) are the quadrature weights of the truncated integral. The
diagonal entries are chosen so that the row sums vanish
\[
(\mathbf M_A)_{ii}
=
-\sum_{j\neq i}(\mathbf M_A)_{ij}.
\]
Since \(\kappa_{ij}=\kappa_{ji}\) and \(\gamma_i+\gamma_j=\gamma_j+\gamma_i\),
the matrix \(\mathbf M_A\) is symmetric. Together with the
symmetry of \(\mathbf L^s\), this implies that the fixed part
\(\mathbf L^s+\mathbf M_A\) is symmetric.

The drift is discretized by a second-order
central difference,
\[
(\mathbf B(\beta))_{ij}
=
b(x_i)\,\frac{\delta_{i+1,j}-\delta_{i-1,j}}{2h},
\]
with one-sided differences at the boundary, and the potential is stored in
the diagonal matrix
\[
\mathbf Q(\theta)=\mathrm{diag}\bigl(q_i+\frac12\sum_j\kappa_{ij}(\gamma_i+\gamma_j)^2\bigr).
\]

Collecting the above, the discrete forward problem reads
\begin{equation}\label{eq:num-system}
\mathbf A(\beta,\theta)\,\mathbf U_R=\mathbf b_R,
\qquad
\mathbf A(\beta,\theta):=
\mathbf L^s+\mathbf M_A
+\mathbf B(\beta)+\mathbf Q(\theta),
\end{equation}
where \(\mathbf b_R\) collects the contribution of the exterior datum \(f\).
We emphasize that \(\mathbf A(\beta,\theta)\) is \emph{affine} in
\((\beta,\theta)\).

\subsection{Observation operator}\label{sec:num-obs}

Let \(W_2\subset\Omega_e\) be the observation set, with points
\(x_\ell\in W_2\), \(\ell=1,\dots,N_o\). Since \(\supp(A)\subseteq\Omega^2\)
and \(W_2\subset\Omega_e\), the magnetic contribution vanishes at the
observation points. The fractional Laplacian at \(x_\ell\) is
\[
(-\Delta)^s u(x_\ell)
=
\sum_{i=1}^{N-1}w_{\ell i}^{\mathrm{int}}u_R(x_i)
+
\sum_{j\in\mathcal E}w_{\ell j}^{\mathrm{ext}}f(x_j),
\]
where \(\mathcal E\) is the index set of the exterior grid points in
\((-R,R)\setminus(-L,L)\). Writing
\(\mathbf L^s_{1}\in\Rb^{N_o\times(N-1)}\) for the matrix with
entries \(w_{\ell i}^{\mathrm{int}}\), and
\(\mathbf b_{e}\in\Rb^{N_o}\) for the vector with entries
\(\sum_{j\in\mathcal E}w_{\ell j}^{\mathrm{ext}}f(x_j)\), the measured
quantity is
\begin{equation}\label{eq:num-observation}
\mathbf h(\beta,\theta)
=
\mathbf L^s_{1}\,\mathbf U_R(\beta,\theta)
+\mathbf b_{e}.
\end{equation}
The magnetic potential affects the data only through its effect on the
interior solution \(\mathbf U_R\).

\subsection{RBF surrogate model}\label{sec:num-rbf}

The repeated evaluation of the forward map \((\beta,\theta)\mapsto\mathbf h\)
during inversion is expensive, since each evaluation requires the solution
of \eqref{eq:num-system}. To decouple the inversion from the forward solve,
we construct an RBF surrogate.

We generate \(N_i\) training parameters
\(\{\boldsymbol\xi^{(m)}\}_{m=1}^{N_i}\subset\Rb^M\), solve \eqref{eq:num-system} for each sample, and collect
\[
\mathbf h^{(m)}=\mathbf h(\boldsymbol\xi^{(m)})\in\Rb^{N_o},
\]
writing \(\mathbf H=[\mathbf h^{(1)},\dots,\mathbf h^{(N_i)}]\in
\Rb^{N_o\times N_i}\).

We use
\[
\varphi_m(\boldsymbol\xi)
=
\frac{1}{\sqrt{\|\boldsymbol\xi-\boldsymbol\xi^{(m)}\|_2^2+\sigma^2}},
\qquad
\sigma
=
\frac{1}{N_i(N_i-1)}
\sum_{i=1}^{N_i}\sum_{j\neq i}
\|\boldsymbol\xi^{(i)}-\boldsymbol\xi^{(j)}\|_2,
\]
and set \(\boldsymbol\varphi(\boldsymbol\xi)
=(\varphi_1(\boldsymbol\xi),\dots,\varphi_{N_i}(\boldsymbol\xi))^\top\).

The coefficient matrix is
\[
\mathbf C=\mathbf H\mathbf F^{-1},
\qquad
\mathbf F_{i,n}=\varphi_n(\boldsymbol\xi^{(i)}).
\]

The surrogate is
\[
\mathbf h_p(\boldsymbol\xi)
=
\mathbf C\boldsymbol\varphi(\boldsymbol\xi)
=
\sum_{n=1}^{N_i}\mathbf c_n\,\varphi_n(\boldsymbol\xi),
\]
where \(\mathbf c_n\) is the \(n\)-th column of \(\mathbf C\). Its Jacobian is
\begin{equation}\label{eq:num-rbf-jac}
\frac{\partial\mathbf h_p}{\partial\xi_j}
=
\mathbf C\,\frac{\partial\boldsymbol\varphi}{\partial\xi_j},
\end{equation}
with
\begin{equation}\label{eq:num-rbf-grad}
\frac{\partial\varphi_m}{\partial\xi_j}
=
-\frac{\xi_j-\xi_j^{(m)}}{\sigma^3(r_m^2+1)^{3/2}},
\qquad
r_m^2=\frac{\|\boldsymbol\xi-\boldsymbol\xi^{(m)}\|_2^2}{\sigma^2}.
\end{equation}

\subsection{Smoothed TV regularization and Levenberg--Marquardt}\label{sec:num-lm}

We estimate \((\beta,\theta)\) by minimizing
\begin{equation}\label{eq:num-objective}
\mathcal J_\varepsilon(\beta,\theta)
=
\frac12\|\mathbf r(\beta,\theta)\|_{L^2}^2
+\alpha_b\sum_{i}\sqrt{(L_b\beta)_i^2+\varepsilon^2}
+\alpha_q\sum_{i}\sqrt{(L_q\theta)_i^2+\varepsilon^2},
\end{equation}
where \(\mathbf r(\beta,\theta):=\mathbf h_p(\beta,\theta)-\mathbf h^\delta\)
is the residual, \((L_b\beta)_i=\beta_{i+1}-\beta_i\), \((L_q\theta)_i=\theta_{i+1}-\theta_i\) are first-difference matrices,
\(\alpha_b,\alpha_q>0\) are regularization parameters, and \(\varepsilon>0\)
is the smoothing parameter. As \(\varepsilon\to0\), the last two terms
converge to the TV semi-norms \(\|L_b\beta\|_1\) and \(\|L_q\theta\|_1\).

We split the objective as
\(\mathcal J_\varepsilon=\Phi+\Psi\), where
\[
\Phi(\beta,\theta)
:=
\frac12\|\mathbf r(\beta,\theta)\|_{L^2}^2,
\qquad
\Psi(\beta,\theta)
:=
\alpha_b\sum_i\sqrt{(L_b\beta)_i^2+\varepsilon^2}
+\alpha_q\sum_i\sqrt{(L_q\theta)_i^2+\varepsilon^2}.
\]
For the data term,
\[
\nabla\Phi(\beta,\theta)
=
\mathbf J^\top\mathbf r,
\qquad
\mathbf J=\frac{\partial\mathbf h_p}{\partial(\beta,\theta)}
\in\Rb^{N_o\times M}.
\]
For the regularization term, set
\(D_b=\mathrm{diag}\bigl(\sqrt{(L_b\beta)_i^2+\varepsilon^2}\bigr)\) and
\(D_q=\mathrm{diag}\bigl(\sqrt{(L_q\theta)_i^2+\varepsilon^2}\bigr)\). Then
\[
\nabla_\beta\Psi=\alpha_b L_b^\top D_b^{-1}L_b\beta,
\qquad
\nabla_\theta\Psi=\alpha_q L_q^\top D_q^{-1}L_q\theta,
\]
so that
\begin{equation}\label{eq:num-gradient}
\nabla\mathcal J_\varepsilon(\beta,\theta)
=
\mathbf J^\top\mathbf r
+
\begin{pmatrix}
\alpha_b L_b^\top D_b^{-1}L_b\beta\\[2pt]
\alpha_q L_q^\top D_q^{-1}L_q\theta
\end{pmatrix}.
\end{equation}

The exact Hessian of the data
term is
\[
\nabla^2 \Phi
=
\mathbf J^\top\mathbf J
+\sum_{\ell=1}^{N_o}
(\mathbf r)_\ell\,\nabla^2(\mathbf h_p)_\ell.
\]
The second term is neglected, as is standard in Gauss--Newton and LM
methods. The Hessian of the regularizer is
\[
\nabla_\beta^2 \Psi
=
\alpha_b L_b^\top
\mathrm{diag}\Bigl(\frac{\varepsilon^2}{((L_b\beta)_i^2+\varepsilon^2)^{3/2}}\Bigr)
L_b,
\]
\[
\nabla_\theta^2 \Psi
=
\alpha_q L_q^\top
\mathrm{diag}\Bigl(\frac{\varepsilon^2}{((L_q\theta)_i^2+\varepsilon^2)^{3/2}}\Bigr)
L_q.
\]
We write
\begin{equation}\label{eq:num-Hreg}
\mathbf H_\varepsilon
:=
\mathrm{diag}\bigl(\nabla_\beta^2 \Psi,\ \nabla_\theta^2 \Psi\bigr),
\end{equation}
so that
\(\nabla^2\mathcal J_\varepsilon\approx
\mathbf J^\top\mathbf J
+\mathbf H_\varepsilon\).

The LM step is obtained by solving the
trust-region subproblem
\[
\min_{\Delta\boldsymbol\xi}\
\nabla\mathcal J_\varepsilon^\top\Delta\boldsymbol\xi
+\frac12\Delta\boldsymbol\xi^\top
\bigl(\mathbf J^\top\mathbf J
+\mathbf H_\varepsilon\bigr)\Delta\boldsymbol\xi
\qquad
\text{s.t.}\quad\|\Delta\boldsymbol\xi\|\le\delta_k.
\]
By the KKT conditions, there exists \(\lambda_k\ge0\) such that
\[
\bigl(\mathbf J^\top\mathbf J
+\mathbf H_\varepsilon
+\lambda_k\mathbf I\bigr)\Delta\boldsymbol\xi
=-\nabla\mathcal J_\varepsilon(\boldsymbol\xi^k),
\qquad
\lambda_k(\|\Delta\boldsymbol\xi\|-\delta_k)=0.
\]
For \(\boldsymbol{\xi}^{k+1} = \boldsymbol{\xi}^k + \Delta\boldsymbol{\xi}\),
the damping parameter \(\lambda_k\) is updated as follows
\begin{equation}\label{eq:num-lambda}
\lambda_{k+1}
=
\begin{cases}
\lambda_k/2, & \text{if }\mathcal J_\varepsilon(\boldsymbol\xi^{k+1})
<\mathcal J_\varepsilon(\boldsymbol\xi^k),\\[2pt]
2\lambda_k, & \text{otherwise.}
\end{cases}
\end{equation}
The iteration is terminated when
\[
\|\nabla\mathcal J_\varepsilon(\boldsymbol\xi^k)\|\le\tau_g,
\qquad\text{or}\qquad
\frac{\|\Delta\boldsymbol\xi\|}{\|\boldsymbol\xi^k\|}\le\tau_x,
\]
with \(\tau_g=\tau_x\) is a constant.

We use the discrepancy principle to select
the regularization parameter. For simplicity we set
\[
\alpha_b=\alpha,\qquad
\alpha_q=\alpha\,\rho,
\qquad
\rho:=\frac{\|L_b\beta_0\|_1}{\|L_q\theta_0\|_1},
\]
where \(\beta_0,\theta_0\) are obtained from a preliminary inversion. The ratio \(\rho\) equalizes the contribution of the two
regularization terms, so that a single parameter \(\alpha\) controls the overall regularization strength. We choose \(\alpha\) as the largest
value such that the weighted residual satisfies
\begin{equation}\label{eq:num-discrepancy}
\|\mathbf h_p(\boldsymbol\xi_\alpha^*)-\mathbf h^\delta\|_{L^2}
\le\tau\delta,
\end{equation}
where \(\tau> 1\) is a constant.

\section{Numerical experiments}\label{sec:num-experiments}

In this section we illustrate the proposed inversion framework with five
numerical examples. All examples are in the known-\(A\) setting, in which
the drift \(b\) and the potential \(q\) can be recovered simultaneously
from exterior DN measurements.

Throughout this section, \(\Omega\) denotes the spatial domain,
\(W_1,W_2\subset\Omega_e\) are non-empty open sets, with \(W_1\) supporting
the exterior data and \(W_2\) the observation set; in the one-dimensional
examples we take \(W_1=W_2\), while in the two-dimensional examples \(W_2\)
is the inner boundary of \(W_1\). To facilitate the calculation, we uniformly assume that \(A_\perp = A_{s\parallel} = 0\). Then, the magnetic flux can be expressed in the following permitted form
\[
A(x,y) = A_{a\parallel}(x,y)
= \alpha(x,y)(\gamma(x) + \gamma(y)).
\]
The profile \(\gamma\) and the order \(s\) are given
in each example.

The drift enters the equation through the term
\(\mathbf b\cdot\nabla u\). In one dimension \(\mathbf b\) reduces to a
scalar \(b(x)\), since \(\nabla u=\partial_x u\) has a single component.
In two dimensions, however, \(\nabla u=(\partial_x u,\partial_y u)\) has
two components, and the drift is necessarily a vector field
\[
\mathbf b(x,y)=\bigl(b_1(x,y),\,b_2(x,y)\bigr),
\]
with \(b_1\) and \(b_2\) the components along the \(x\) and \(y\) axes.
The drift term is then the inner product
\[
\mathbf b\cdot\nabla u
=
b_1\,\partial_x u+b_2\,\partial_y u.
\]
Both components are recovered simultaneously in the two-dimensional
examples below.

Since \(\supp(A)\subseteq\Omega^2\), the magnetic contribution vanishes on
\(W_2\), and the observation operator is
\[
\mathcal O(u)=(-\Delta)^s u\big|_{W_2}.
\]
The measurements are all affected by additive Gaussian noise with standard deviation \(\delta\).

\begin{figure}[htbp]
	\centering
    \captionsetup{font={small}}
	{
		\begin{minipage}{0.45\linewidth}
			\centering
			\includegraphics[height=0.7\columnwidth,width=0.9\linewidth]{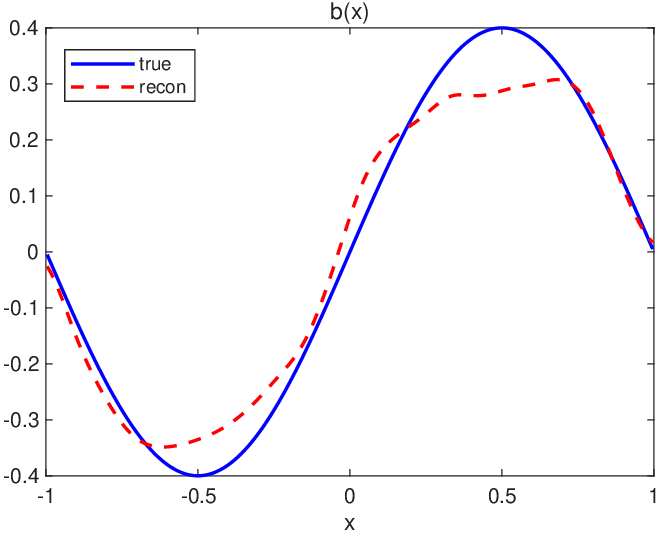}\\
			(a)
		\end{minipage}
		\begin{minipage}{0.45\linewidth}
			\centering
			\includegraphics[height=0.7\columnwidth,width=0.9\linewidth]{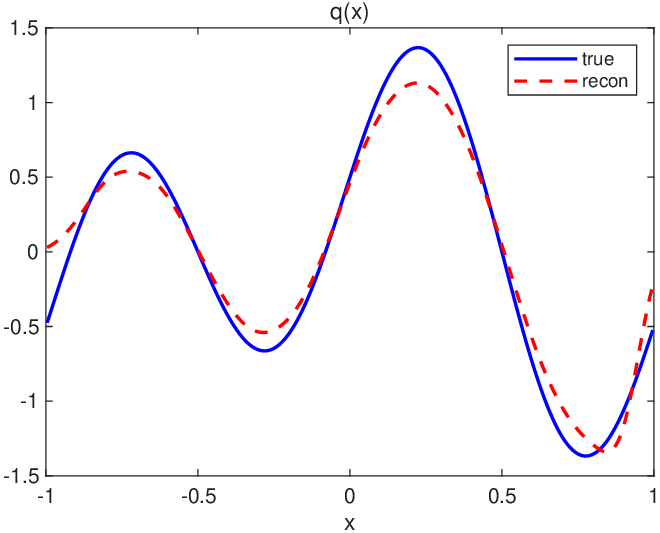}\\
			(b)
		\end{minipage}
	}
\caption{The numerical results for Example 1.}
	\label{fig1}
\end{figure}

\emph{Example 1.} We consider that both $b$ and $q$ are smooth functions,
\[
b(x)=0.4\sin(\pi x),
\qquad
q(x)=\sin(2\pi x)+0.5\cos(\pi x).
\]
The magnetic potential is fixed and given by
\begin{equation}\label{eq:ex1-A}
A(x,y)=\alpha(x,y)\bigl(\gamma(x)+\gamma(y)\bigr),
\qquad
\gamma(x)=0.3\sin(\pi x),
\end{equation}
where
\[
\alpha(x,y)
=
\frac{c_{1,s}^{1/2}}{\sqrt2}
\frac{y-x}{|y-x|^{1/2+s+1}},
\qquad s=0.7.
\]
In Figure \ref{fig1}, Figures (a) and (b) respectively present the inversion images of the drift term $b$ and the potential term $q$ when the error level is $10^{-4}$.

\begin{figure}[htbp]
	\centering
    \captionsetup{font={small}}
	{
		\begin{minipage}{0.45\linewidth}
			\centering
			\includegraphics[height=0.7\columnwidth,width=0.9\linewidth]{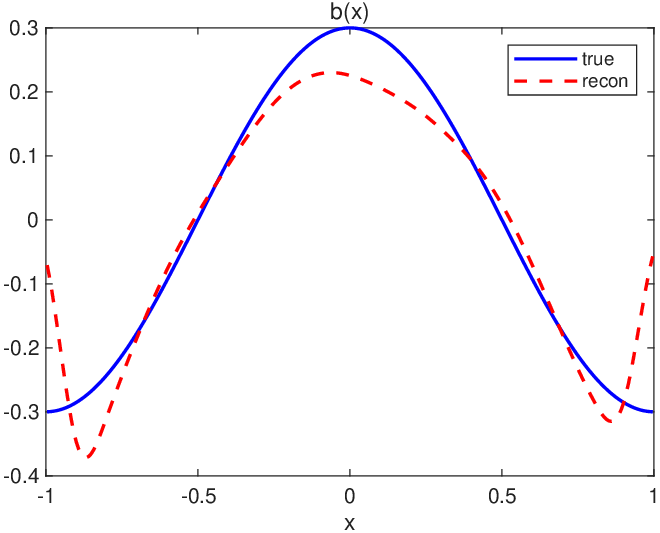}\\
			(a)
		\end{minipage}
		\begin{minipage}{0.45\linewidth}
			\centering
			\includegraphics[height=0.7\columnwidth,width=0.9\linewidth]{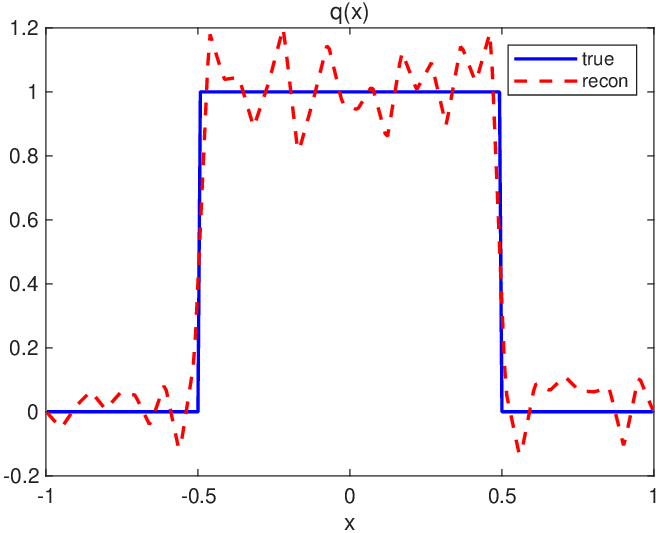}\\
			(b)
		\end{minipage}
	}
\caption{The numerical results for Example 2.}
	\label{fig2}
\end{figure}

\emph{Example 2.} We consider the following example and take
\(s=0.8\),
\[
b(x)=0.3\cos(\pi x),
\]
and the piecewise constant potential
\[
q(x)=
\begin{cases}
0, & x\in(-1,-0.5),\\
1, & x\in(-0.5,0.5),\\
0, & x\in(0.5,1).
\end{cases}
\]
The magnetic potential is fixed and given by
\begin{equation}\label{eq:ex2-A}
A(x,y)=\alpha(x,y)\bigl(\gamma(x)+\gamma(y)\bigr),
\qquad
\gamma(x)=0.2\cos(\pi x)(1-x^2).
\end{equation}
 Figure~\ref{fig2} shows that the inversion images of the drift term $b$ and the potential term $q$ when the error level is $10^{-4}$.

\begin{figure}[htbp]
	\centering
    \captionsetup{font={small}}
	{
		\begin{minipage}{0.45\linewidth}
			\centering
			\includegraphics[height=0.7\columnwidth,width=0.9\linewidth]{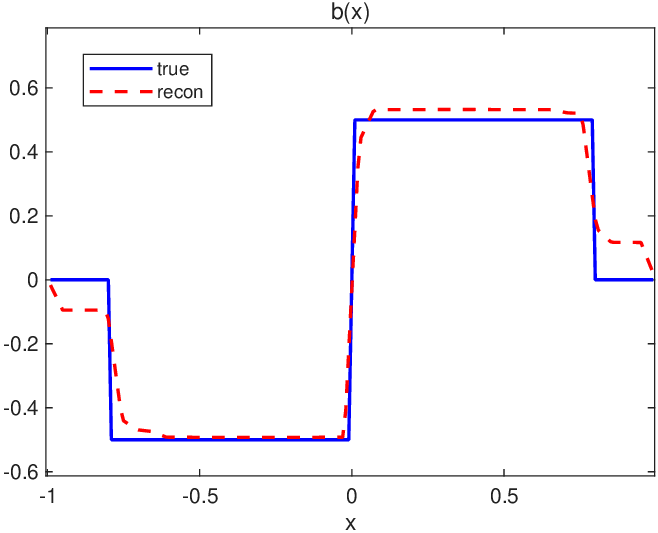}\\
			(a)
		\end{minipage}
		\begin{minipage}{0.45\linewidth}
			\centering
			\includegraphics[height=0.7\columnwidth,width=0.9\linewidth]{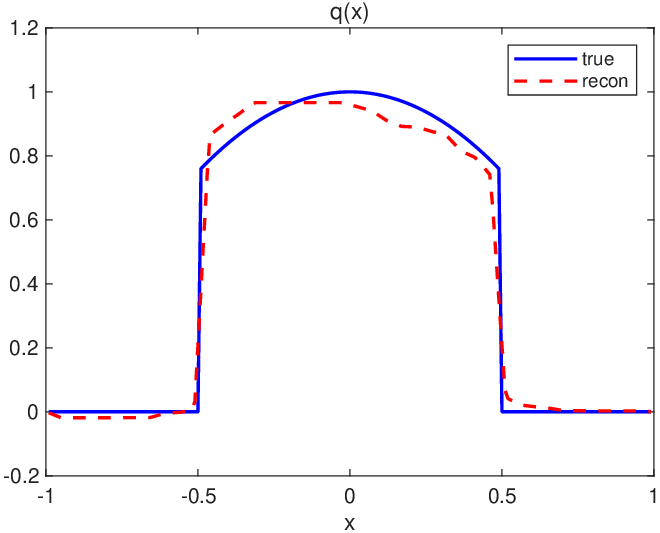}\\
			(b)
		\end{minipage}
	}
\caption{The numerical results for Example 3.}
	\label{fig3}
\end{figure}

\emph{Example 3.}We consider that both $b$ and $q$ are discontinuous functions,
\[
b(x)=
\begin{cases}
-0.5, & x\in(-0.8,0),\\
0.5, & x\in(0,0.8),\\
0, & x\in(-1,-0.8]\cup[0.8,1),
\end{cases}
\]
\[
q(x)=
\begin{cases}
1-x^2, & x\in(-0.5,0.5),\\
0, & x\in(-1,-0.5]\cup[0.5,1).
\end{cases}
\]
In this example, let $s=0.8$, the magnetic potential is given by $\gamma(x)=0.2\sin(2\pi x)\,(1-x^2)$. Figure~\ref{fig3} shows that the inversion images of the drift term $b$ and the potential term $q$ when the error level is $\delta=10^{-5}$. This example confirms that the smoothed TV regularization handles multiple jumps simultaneously in both coefficients.

\begin{figure}[htbp]
	\centering
    \captionsetup{font={small}}
	{
		\begin{minipage}{0.3\linewidth}
			\centering
			\includegraphics[height=0.7\columnwidth,width=0.9\linewidth]{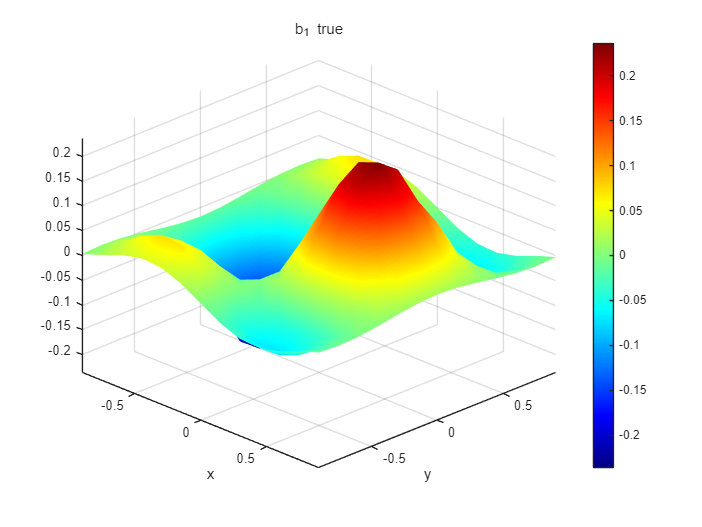}\\
		\end{minipage}
		\begin{minipage}{0.3\linewidth}
			\centering
			\includegraphics[height=0.7\columnwidth,width=0.9\linewidth]{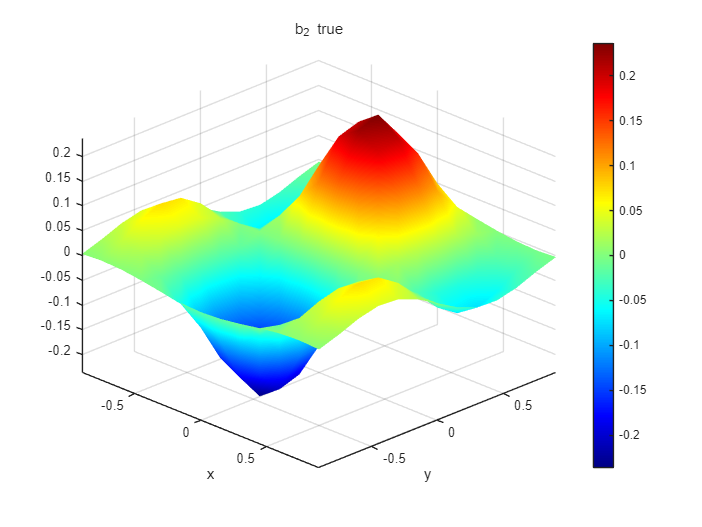}\\
		\end{minipage}
		\begin{minipage}{0.3\linewidth}
			\centering
			\includegraphics[height=0.7\columnwidth,width=0.9\linewidth]{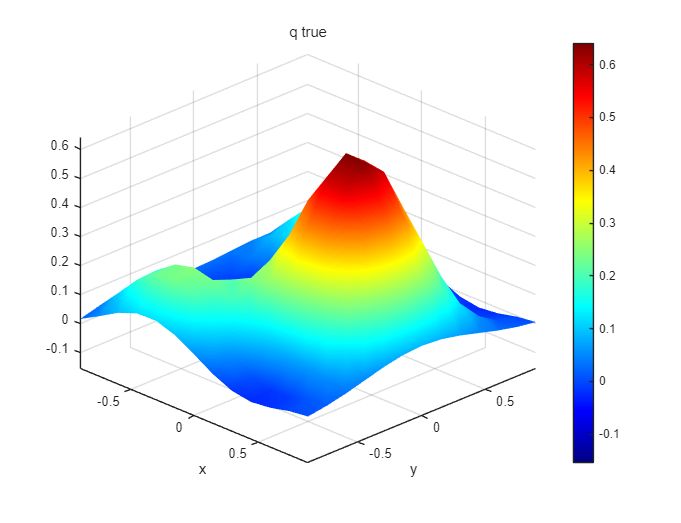}\\
		\end{minipage}
	}
\caption{The real solution of Example 4.}
	\label{fig41}
\end{figure}

\begin{figure}[htbp]
	\centering
    \captionsetup{font={small}}
	{
		\begin{minipage}{0.3\linewidth}
			\centering
			\includegraphics[height=0.7\columnwidth,width=0.9\linewidth]{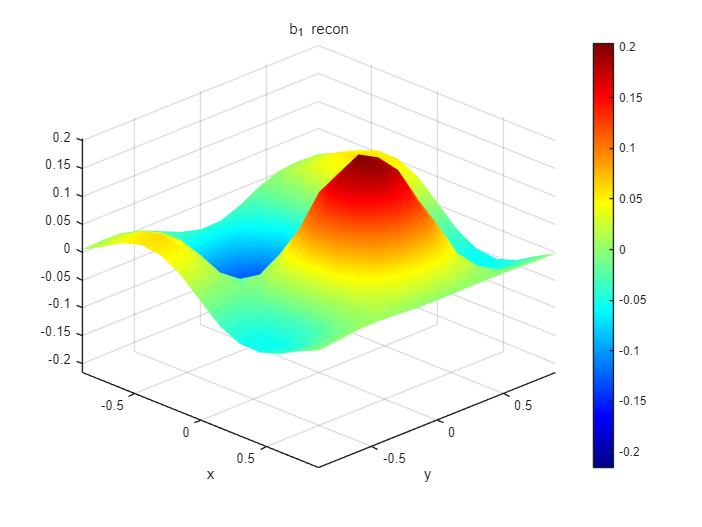}\\
		\end{minipage}
		\begin{minipage}{0.3\linewidth}
			\centering
			\includegraphics[height=0.7\columnwidth,width=0.9\linewidth]{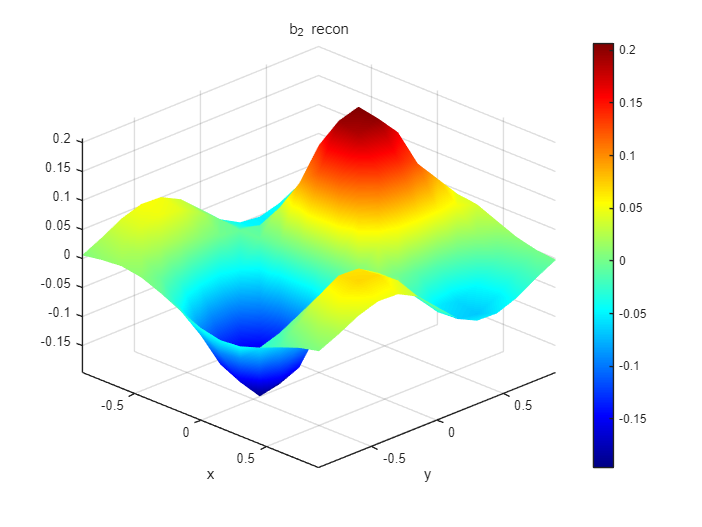}\\
		\end{minipage}
		\begin{minipage}{0.3\linewidth}
			\centering
			\includegraphics[height=0.7\columnwidth,width=0.9\linewidth]{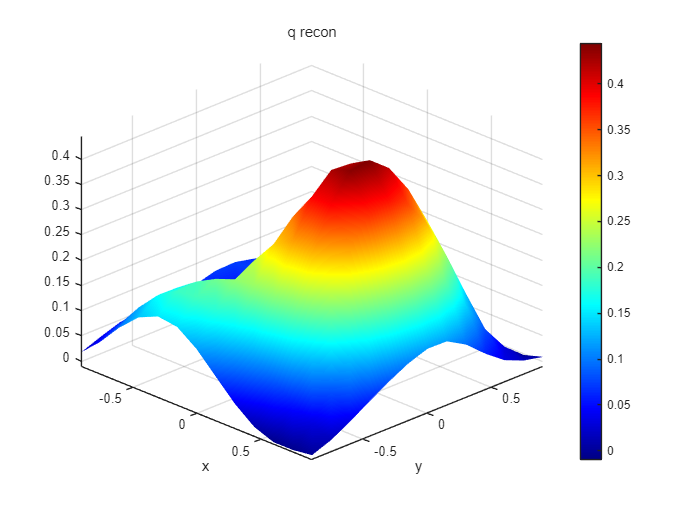}\\
		\end{minipage}
	}
\caption{The reconstructed solution of Example 4 at the error level of \(\delta=10^{-5}\).}
	\label{fig42}
\end{figure}

\begin{figure}[htbp]
	\centering
    \captionsetup{font={small}}
	{
		\begin{minipage}{0.3\linewidth}
			\centering
			\includegraphics[height=0.7\columnwidth,width=0.9\linewidth]{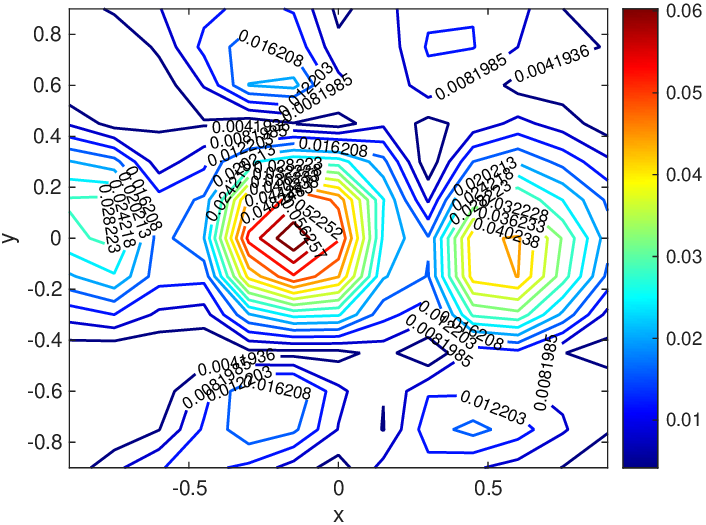}\\
		\end{minipage}
		\begin{minipage}{0.3\linewidth}
			\centering
			\includegraphics[height=0.7\columnwidth,width=0.9\linewidth]{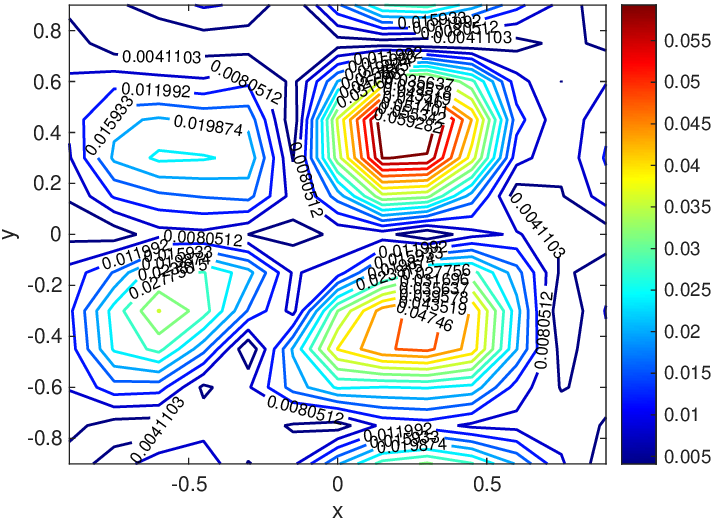}\\
		\end{minipage}
		\begin{minipage}{0.3\linewidth}
			\centering
			\includegraphics[height=0.7\columnwidth,width=0.9\linewidth]{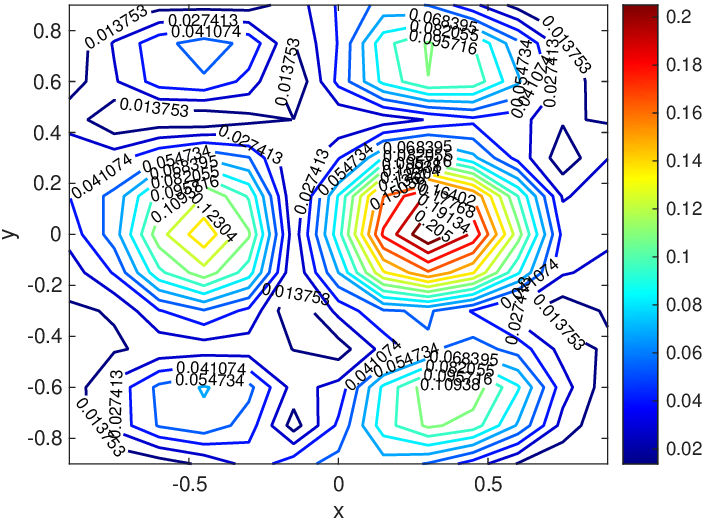}\\
		\end{minipage}
	}
\caption{Contour plot of the relative error of Example 4.}
	\label{fig43}
\end{figure}

\emph{Example 4.} We consider the following two-dimensional numerical example,
\[
b_1(x,y)=0.3\sin(\pi x)\cos(\pi y)(1-x^2)(1-y^2),\quad
b_2(x,y)=0.3\cos(\pi x)\sin(\pi y)(1-x^2)(1-y^2),
\]
\[
q(x,y)=0.5\sin(\pi x)\cos(\pi y)(1-x^2)(1-y^2)+0.3(1-x^2)(1-y^2).
\]
Let \(\Omega=(-1,1)^2\), \(s=0.6\), \(R=3\). The exterior data is \(W_1=(-3,3)^2\setminus(-1-\epsilon,1+\epsilon)^2\), \(\epsilon=0.05\),
and the observation set is the inner boundary
\[
W_2=\{(x,y)\in W_1:\ |x|=1+\epsilon\text{ or }|y|=1+\epsilon\}.
\]
Figure~\ref{fig41} shows the exact solutions of $b_1$, $b_2$, and $q$. Figure~\ref{fig42} presents the reconstructed solutions of  $b_1$, $b_2$, and $q$ when the error is $\delta=10^{-5}$. The relative error computed for the exact solution against its reconstruction is illustrated in Figure \ref{fig43}.


\begin{figure}[htbp]
	\centering
    \captionsetup{font={small}}
	{
		\begin{minipage}{0.3\linewidth}
			\centering
			\includegraphics[height=0.7\columnwidth,width=0.9\linewidth]{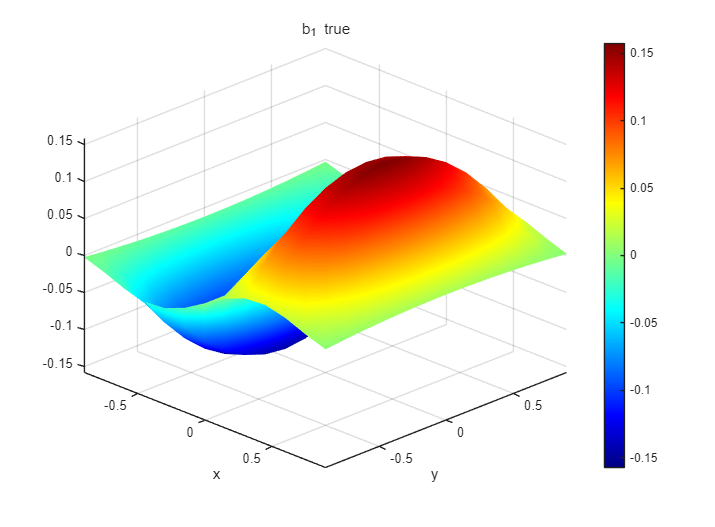}\\
		\end{minipage}
		\begin{minipage}{0.3\linewidth}
			\centering
			\includegraphics[height=0.7\columnwidth,width=0.9\linewidth]{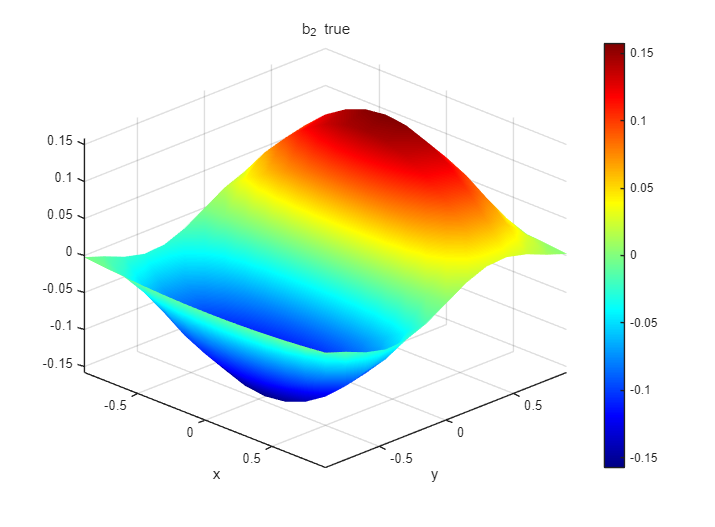}\\
		\end{minipage}
		\begin{minipage}{0.3\linewidth}
			\centering
			\includegraphics[height=0.7\columnwidth,width=0.9\linewidth]{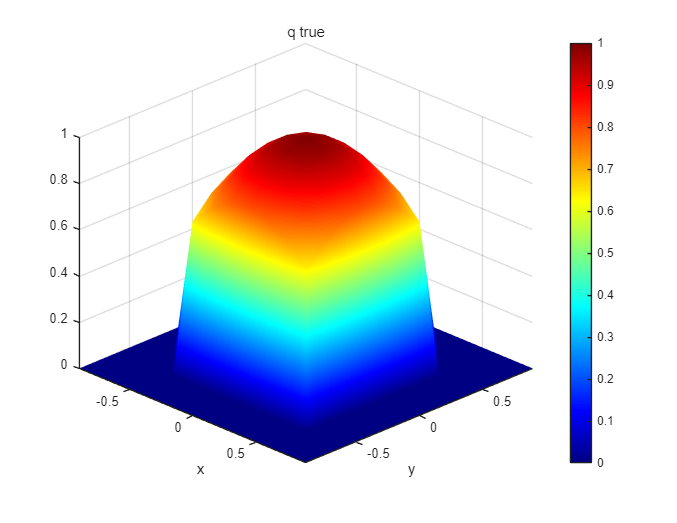}\\
		\end{minipage}
	}
\caption{The real solution of Example 5.}
	\label{fig51}
\end{figure}

\begin{figure}[htbp]
	\centering
    \captionsetup{font={small}}
	{
		\begin{minipage}{0.3\linewidth}
			\centering
			\includegraphics[height=0.7\columnwidth,width=0.9\linewidth]{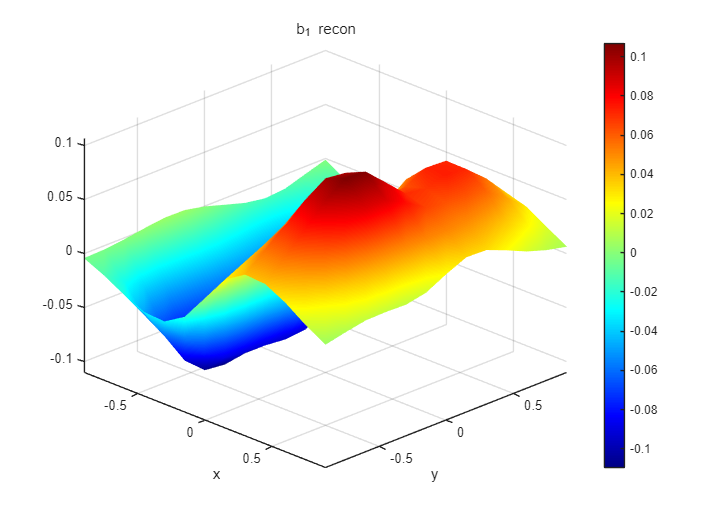}\\
		\end{minipage}
		\begin{minipage}{0.3\linewidth}
			\centering
			\includegraphics[height=0.7\columnwidth,width=0.9\linewidth]{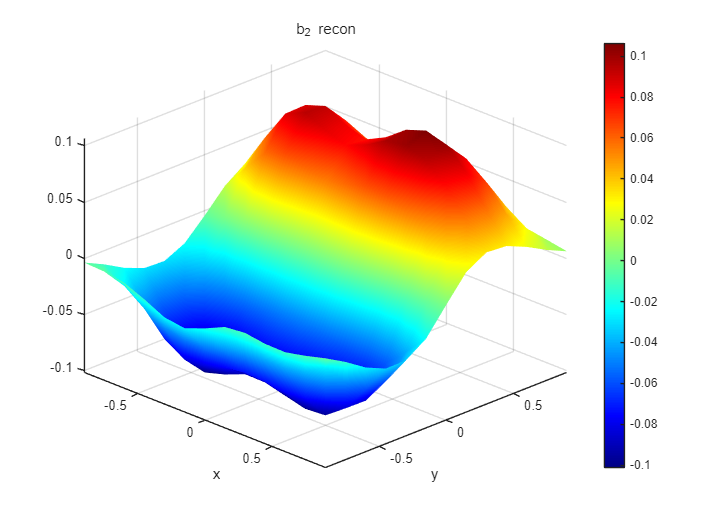}\\
		\end{minipage}
		\begin{minipage}{0.3\linewidth}
			\centering
			\includegraphics[height=0.7\columnwidth,width=0.9\linewidth]{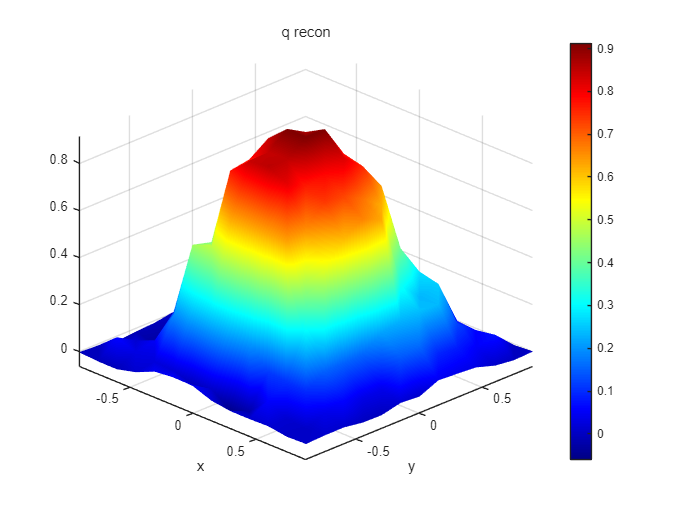}\\
		\end{minipage}
	}
\caption{The reconstructed solution of Example 5 at the error level of \(\delta=10^{-4}\).}
	\label{fig52}
\end{figure}

\begin{figure}[htbp]
	\centering
    \captionsetup{font={small}}
	{
		\begin{minipage}{0.3\linewidth}
			\centering
			\includegraphics[height=0.7\columnwidth,width=0.9\linewidth]{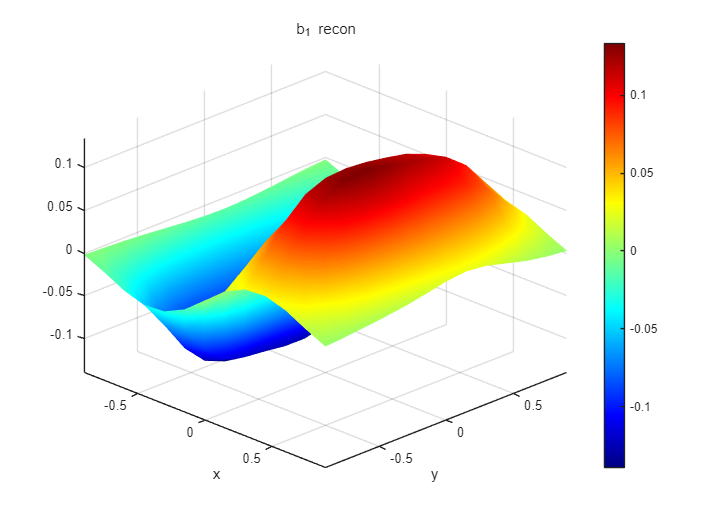}\\
		\end{minipage}
		\begin{minipage}{0.3\linewidth}
			\centering
			\includegraphics[height=0.7\columnwidth,width=0.9\linewidth]{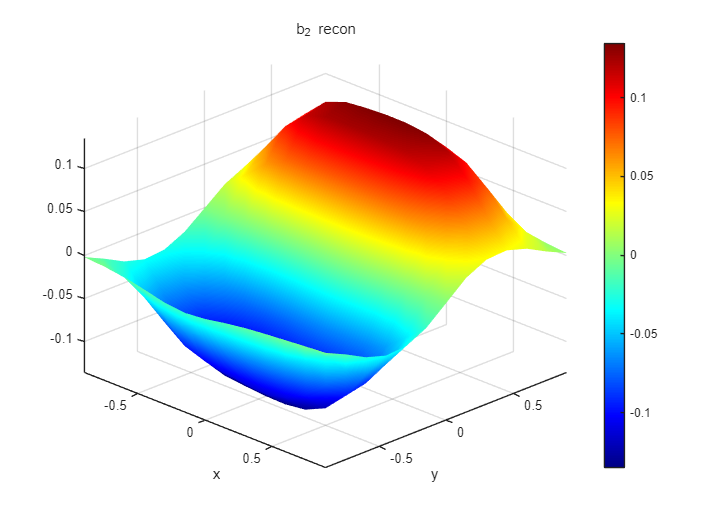}\\
		\end{minipage}
		\begin{minipage}{0.3\linewidth}
			\centering
			\includegraphics[height=0.7\columnwidth,width=0.9\linewidth]{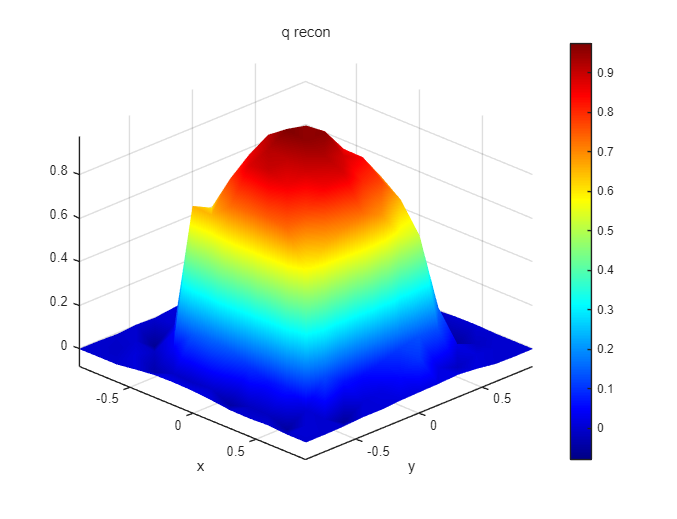}\\
		\end{minipage}
	}
\caption{The reconstructed solution of Example 5 at the error level of \(\delta=10^{-5}\).}
	\label{fig53}
\end{figure}

\begin{figure}[htbp]
	\centering
    \captionsetup{font={small}}
	{
		\begin{minipage}{0.3\linewidth}
			\centering
			\includegraphics[height=0.7\columnwidth,width=0.9\linewidth]{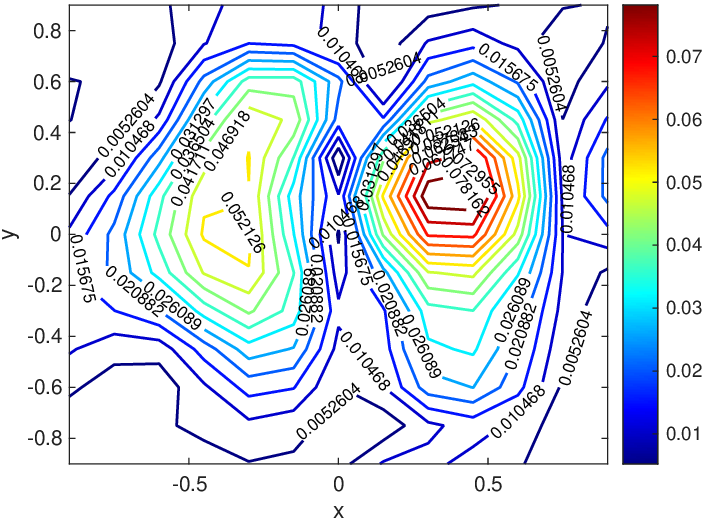}\\
		\end{minipage}
		\begin{minipage}{0.3\linewidth}
			\centering
			\includegraphics[height=0.7\columnwidth,width=0.9\linewidth]{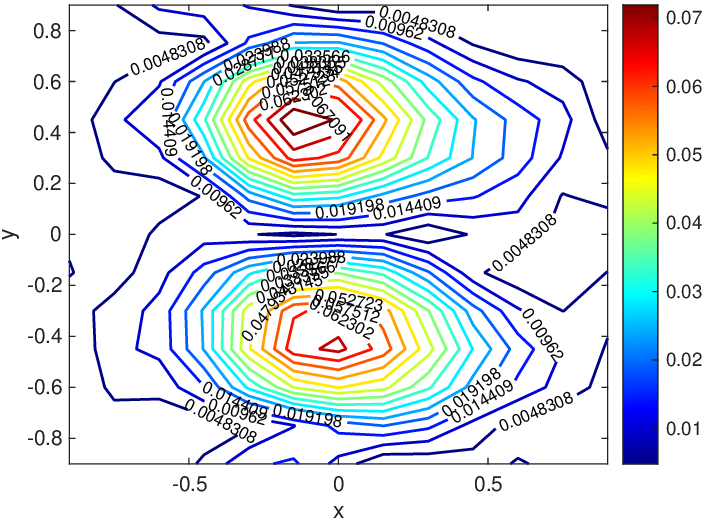}\\
		\end{minipage}
		\begin{minipage}{0.3\linewidth}
			\centering
			\includegraphics[height=0.7\columnwidth,width=0.9\linewidth]{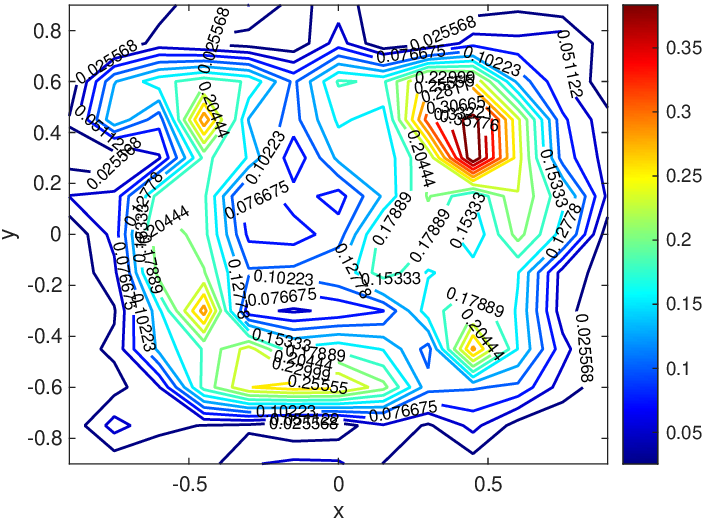}\\
		\end{minipage}
	}
\caption{Contour plot of the relative error of Example 5 at \(\delta=10^{-4}\).}
	\label{fig54}
\end{figure}

\begin{figure}[htbp]
	\centering
    \captionsetup{font={small}}
	{
		\begin{minipage}{0.3\linewidth}
			\centering
			\includegraphics[height=0.7\columnwidth,width=0.9\linewidth]{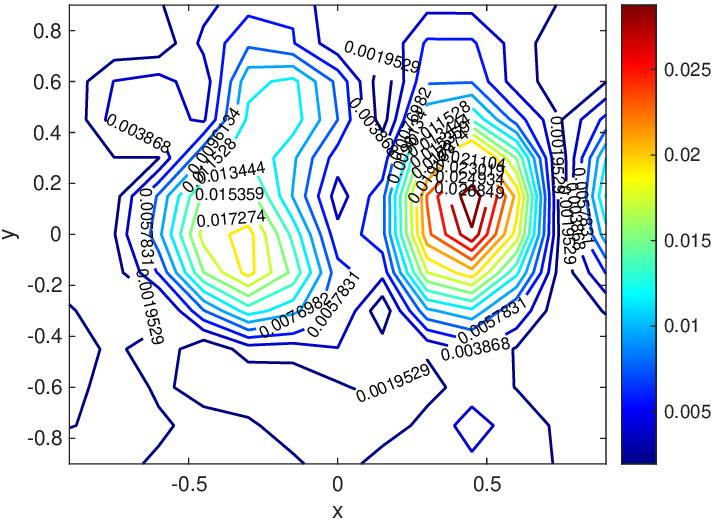}\\
		\end{minipage}
		\begin{minipage}{0.3\linewidth}
			\centering
			\includegraphics[height=0.7\columnwidth,width=0.9\linewidth]{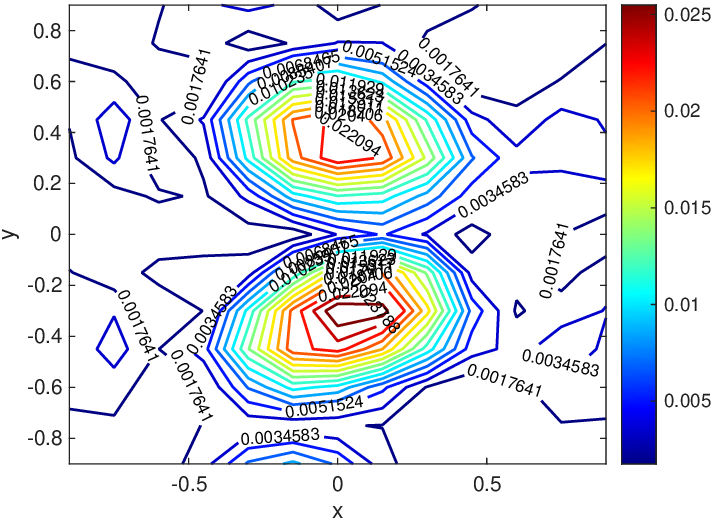}\\
		\end{minipage}
		\begin{minipage}{0.3\linewidth}
			\centering
			\includegraphics[height=0.7\columnwidth,width=0.9\linewidth]{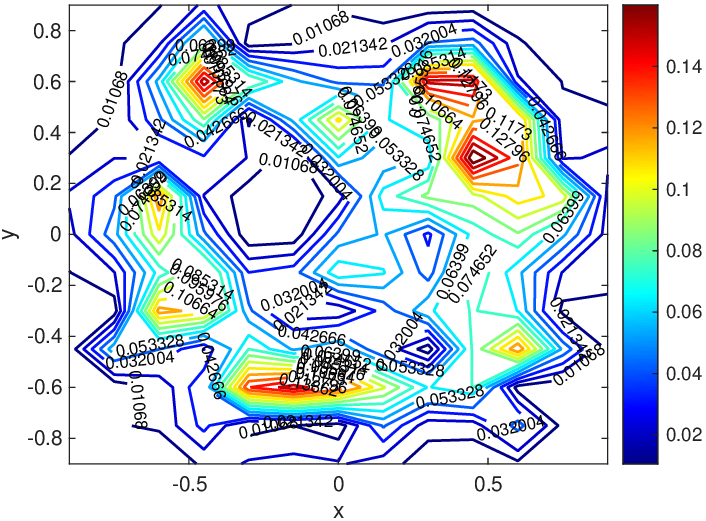}\\
		\end{minipage}
	}
\caption{Contour plot of the relative error of Example 5 at \(\delta=10^{-5}\).}
	\label{fig55}
\end{figure}

\emph{Example 5.} Two dimensions, discontinuous potential.
\[
b_1(x,y)=0.2\sin(\pi x)(1-x^2)(1-y^2),\quad
b_2(x,y)=0.2\sin(\pi y)(1-x^2)(1-y^2),
\]
\[
q(x,y)=
\begin{cases}
(1-x^2)(1-y^2), & (x,y)\in(-0.5,0.5)^2,\\
0, & \text{otherwise},
\end{cases}
\]
 We use the same domain, exterior data, observation set as in Example~4, and let \(s=0.7\). Figure~\ref{fig51} shows the exact solutions of $b_1$, $b_2$, and $q$. Figures ~\ref{fig52} and ~\ref{fig53} respectively illustrate the reconstructed solutions of $b_1$, $b_2$ and $q$ when the error is $\delta = 10^{-4}$ and $\delta = 10^{-5}$. At the above two error levels, the relative error computed for the exact solution against its reconstruction is illustrated in Figures ~\ref{fig54} and ~\ref{fig55} respectively. Meanwhile, Table \ref{tab51} presents the $L^2$-error of $(b1, b2, q)$ under different orders $s$ and error levels $\delta$.

\begin{table}	
	\centering
	\tabcolsep=0.1cm	\caption{The $L^2$ error of ($b_1$, $b_2$, $q$) under different fractional order $s$ and error level $\delta$ in Example 5.}
	\begin{tabular}{c|ccc}	
		\hline
		$s$  $\backslash$ $\delta$  & $\delta=10^{-4}$ & $\delta=10^{-5}$ & $\delta=10^{-6}$\\
		\hline
		$s=0.3$	& $(0.1980,0.1791,0.1666)$ & $(0.0712,0.0591,0.0720)$  & $(0.0252,0.0253,0.0522)$ \\
		$s=0.5$	& $(0.2784,0.2516,0.2359)$ & $(0.0740,0.0541,0.0772)$  & $(0.0296,0.0307,0.0627)$ \\
		$s=0.7$	& $(0.3448,0.3189,0.2917)$ & $(0.0842,0.0559,0.1952)$  & $(0.0332,0.0358,0.0595)$ \\
		\hline
	\end{tabular}
	\label{tab51}
\end{table}


As illustrated by the numerical examples above, our inversion strategy performs satisfactorily. The RBF surrogate model technique brings about a considerable reduction in computational expense and thus enhances the computational efficiency dramatically. Meanwhile, the outcomes confirm the validity of the uniqueness theory established for the simultaneous inversion of the drift magnetic fractional Schr\"odinger equation.
\section{Conclusion}\label{sec:concl}

In this paper we have studied the Calder\'on problem for the fractional
magnetic Schr\"odinger equation with drift,
\[
  (-\Delta)_A^s u+b\cdot\nabla u+q u=0 \quad\text{in }\Omega,
\]
in the bivariate magnetic potential framework.
We proved that if the magnetic potential \(A\) is known and fixed, then the drift \(b\) and the scalar potential \(q\) are simultaneously and uniquely
determined by exterior measurements. If \(A\) is unknown, then the triple \((A,b,q)\) can only be
recovered  by a gauge equivalence relation \(\sim\). The proofs rely on the Alessandrini identity, the weak unique continuation
property, and the Runge approximation property. Meanwhile, we have developed a deterministic inversion
framework for the known-\(A\) setting. The forward problem is discretized by a finite-difference scheme on a truncated domain, the unknown coefficients
\(b\) and \(q\) are parameterized, and an RBF surrogate is constructed. The inversion is performed by a Levenberg--Marquardt method applied to a smoothed Tikhonov--TV objective,
with the regularization parameter selected by the discrepancy principle. Numerical experiments in one and two dimensions
confirm the theoretical uniqueness result.

\bibliographystyle{unsrt}
\bibliography{liu}

\end{document}